\documentclass{article}
\usepackage{arxiv}

\usepackage{subfig}
\usepackage{lineno,hyperref}
\modulolinenumbers[5]

\usepackage{multirow}
\usepackage{textcomp}
\usepackage{stmaryrd}
\usepackage{boldfonts}
\usepackage{symbols}
\usepackage{ulem}
\usepackage{lipsum}
\usepackage{amsfonts}
\usepackage{graphics}
\usepackage{graphicx}
\usepackage{epstopdf}

\usepackage[ruled,vlined]{algorithm2e}
\usepackage{algorithmic}

\usepackage{diagbox}
\usepackage{mathtools}
\usepackage{tikz}
\ifpdf
  \DeclareGraphicsExtensions{.eps,.pdf,.png,.jpg}
\else
  \DeclareGraphicsExtensions{.eps}
\fi

\usepackage{amsopn}

\usepackage{eqnarray}

\newcommand{\vertiii}[1]{{\left\vert\kern-0.25ex\left\vert\kern-0.25ex\left\vert #1 
    \right\vert\kern-0.25ex\right\vert\kern-0.25ex\right\vert}}

\title{Choice of Interior Penalty Coefficient for Interior Penalty Discontinuous Galerkin Method for Biot's System by Employing Machine Learning}

\author{ Sanghyun Lee \\ 
Department of Mathematics \\ Florida State University \\ United States of America \\
  \texttt{lee@math.fsu.edu} 
 \And
Teeratorn Kadeethum \\
The Danish Hydrocarbon Research and Technology Centre \\ 
Technical University of Denmark \\ 
 Denmark 
\And 
Hamidreza M. Nick\\
The Danish Hydrocarbon Research and Technology Centre \\ 
Technical University of Denmark \\ 
 Denmark 
}

\begin{document}
\maketitle

\begin{abstract}
  In this paper, the optimal choice of the interior penalty parameter of the discontinuous Galerkin finite element methods for both the elliptic problems and the Biot's systems are studied by utilizing the neural network and machine learning. It is crucial to choose the optimal interior penalty parameter, which is not too small or not too large for the stability, robustness, and efficiency of the numerical discretized solutions. Both linear regression and nonlinear artificial neural network methods are employed and compared using several numerical experiments to illustrate the capability of our proposed computational framework. This framework is an integral part of a developing automated numerical simulation platform because it can automatically identify the optimal interior penalty parameter. Real-time feedback could also be implemented to update and improve model accuracy on the fly. 
\end{abstract}


\section{Introduction}

Discontinuous Galerkin finite element method (DG) is one of the most popular non conforming finite elements employed for various realistic applications, especially with discontinuous coefficients. The idea of DG finite element methods originated from \cite{nitsche1971variationsprinzip} and extended by several authors, including \cite{douglas1976interior,wheeler1978elliptic,percell1978local,arnold1982interior}, which were also called Interior Penalty Galerkin Methods.  
DG has been actively employed in many multiphysics applications due to the following advantages. 
First, DG is one of the well known and successful methods in terms of the local flux conservation with highly varying material properties \cite{ref1,riviere2000discontinuous,CNM:CNM464,cockburn1999some,cockburn1998local}. 
In addition, DG can deal robustly with general partial differential equations as well as with equations whose type changes within the computational domain, such as from advection dominated to diffusion dominated \cite{Sun2005,Sun2006,babuvska1999discontinuous}.

However, one of the main disadvantages of DG is that the stability and the accuracy of the scheme depend on the interior penalty parameter that needs to be chosen. Numerical analyses of DG are proved under an assumption on the interior penalty parameter, and it is crucial to employ the optimal interior penalty parameter. 
Generally, if the parameter is too large, DG schemes converge to the continuous Galerkin finite element methods and often suffer from the linear solver. If the parameter is too small, the stability of the scheme is not guaranteed. 
Thus, several studies of the lower bounds for the penalty parameter have been obtained in \cite{ainsworth2007posteriori,ainsworth2010fully,ainsworth2009constant,epshteyn2007estimation,shahbazi2005explicit}. 
Moreover, weighted interior penalty parameters for the cases where the diffusion coefficient is discontinuous were studied in \cite{ern2008posteriori,ern2009discontinuous}, and specific illustrations on the selection of the penalty parameters are shown in \cite{ainsworth2012note}.

In this paper, we propose a new procedure to find the optimal interior penalty parameters for both elliptic problems and the poroelastic Biot system. Since the choice of the optimal interior penalty parameters for multiphysics multiscale coupled problems or problems with discontinuous and heterogeneous material properties are nontrivial by the traditional analytic approaches, we employ machine learning processes to predict the optimal interior penalty parameters. 
Many machine learning models have been a center of attention for the past decades because of its approximation power that could be practically applied to various applications~\cite{libbrecht2015machine,brink2016real}.
These algorithms range from classic linear regression models \cite{thomas1999jmp,seber2012linear}, spatial interpolation techniques such as kriging \cite{cressie1988spatial} or maximum likelihood estimation \cite{myung2003tutorial}, and nonlinear approximation functions such nonlinear regression \cite{park2001nonlinear} or deep learning \cite{goodfellow2016deep}. Note that these methods are only used as examples and by no means are the completed set of available algorithms. 

Recently, deep learning has become more attractive because it is scalable \cite{chilimbi2014project}, suitable for GPU functionality \cite{cui2016geeps}, and required computational resources become less demanding because of the mini-batch gradient descent approach \cite{hinton2012neural}. Deep learning has also been successfully applied to solve partial differential equations, which generally are solved by classical numerical methods such as finite difference, finite volume, or finite element methods \cite{wang2017physics,raissi2019physics}. Moreover, this technique has been used to assist the traditional numerical methods such as finite element to enhance its performance \cite{oishi2017computational}. Hence, this paper aims to apply this method for identifying the optimal interior penalty parameters in complex problems.

The proposed procedure benefits not only the simple elliptic problem or Biot's equations but also any multiphysics multiscale coupled problems. Besides, in cases where many simulations have to be performed with different settings, e.g., mesh size, material properties, or various interior penalty schemes, our proposed framework can automatically identify the optimal interior penalty parameter. Real-time feedback could also be implemented to update and improve model accuracy. 

The paper is organized as follows. Our governing system and finite element discretizations are in
Section \ref{sec:main}. Details about the machine learning algorithm are discussed in Section  \ref{sec:alg}.  The numerical results are in Section  \ref{sec:experiments}; this section illustrates the effects of interior penalty parameters on both solution quality and simulation behavior. Performance between linear and nonlinear approximation functions are also compared. Finally, the conclusions follow in Section  \ref{sec:conclusions}.

\section{Mathematical Model}
\label{sec:main}

In this section, we briefly recapitulate the Biot system for poro-elasticity that we will discuss in this paper.  
Let $\Omega \subset \mathbb{R}^d$ ($d \in \{1,2,3\}$) be the computational domain, which is bounded by the boundary, $\partial \Omega$. The time domain is denoted by $\mathbb{T} = \left(0,T\right]$ with $T>0$. 
Then the coupling between the fluid flow and solid deformation can be captured through the application of Biot\textquotesingle s equation of poroelasticity, which is composed of linear momentum and mass balance equations \cite{biot1941general}.

First, the mass balance equation is given as \cite{coussy2004poromechanics}:
\begin{equation} \label{eq:mass_balance}
\rho\left(\phi c_{f}+\dfrac{\alpha-\phi}{K_{s}}\right) \dfrac{\partial}{\partial t}p + 
\rho \alpha \frac{\partial}{\partial t} \nabla \cdot \bu -
\nabla \cdot \bkappa(\nabla p-\rho \mathbf{g})=g \text { in } \Omega \times \mathbb{T},
\end{equation}
where 
$p(\cdot ,  t) : \Omega \times  (0; T] \to \mathbb{R}$  
is a scalar-valued fluid pressure, 
$\mathbf{u} (\cdot , t) : \Omega \times  (0; T] \to \mathbb{R}^d$ 
is a vector-valued displacement,
$\rho$ is a fluid density, $\phi$ is an initial porosity, 
$c_f$ is a fluid compressibility, $\bg$ is a gravitational vector, $g$ is a sink/source.
Here, $\nabla \cdot \bu$ term represents the volumetric deformation and $\bkappa$ is defined as:
\begin{equation}
\bkappa:=\frac{\rho {\bk_{m}}}{\mu}, 
\label{eq:kappa}
\end{equation}
where $\bk_m$ is a matrix permeability tensor and $\mu$ is a fluid viscosity.




The mass balance equation (the fluid flow problem) is supplemented by the following boundary and initial conditions: 
\begin{eqnarray}
p&=&p_{D} \text { on } \partial \Omega_{p} \times \mathbb{T}, \\
-\nabla \cdot \bkappa(\nabla p-\rho \mathbf{g}) \cdot \bn&=&q_{D} \text { on } \partial \Omega_{q} \times \mathbb{T}, \\
p&=&p_{0} \text { in } \Omega \text { at } t = 0,
\end{eqnarray}
where $p_D$ and $q_D$ are specified pressure and flux, respectively,
and $\partial \Omega$ is decomposed to pressure and flux boundaries, $\partial \Omega_p$ and $\partial \Omega_q$, respectively.

Secondly, the linear momentum balance equation can be written as follows: 
\begin{equation}
\nabla \cdot \bsigma (\bu,p) =\bf{f}.
\end{equation}
For the simplicity,  a body force $\bf{f}$ is neglected in this study.  Here, $\bsigma$ is total stress, which is defined as:
\begin{equation}
\bsigma := \bsigma(\bu, p) = \bsigma^{\prime}(\bu) - \alpha p \bI,
\end{equation}
where $\bI$ is the identity tensor and $\alpha$ is Biot\textquotesingle s coefficient defined as \cite{Jaeger2010}:
\begin{equation} \label{eq:biot_coeff}
\alpha:=1-\frac{K}{K_{{s}}},
\end{equation}
with the bulk modulus of a rock matrix $K$ and the solid grains modulus $K_s$. In addition, $\bsigma^{\prime}$ is an effective stress written as: 
\begin{equation}
\bsigma^{\prime}:= \bsigma^{\prime}(\bu) = 
2 \mu_{l} \bepsilon(\bu)-\lambda_{l} \nabla \cdot \bu \bI,
\end{equation}
where $\lambda_{l}$ and $\mu_{l}$ are Lam\'e constants. By assuming a small displacement, a strain is defined as: 

\begin{equation}
\bepsilon(\bu) :=\frac{1}{2}\left(\nabla \bu+\nabla \bu^{T}\right),
\end{equation}

Thus, we can write the linear momentum balance supplemented by its boundary and initial conditions as:
\begin{eqnarray}
\label{eq:linear_balance}
\nabla \cdot \bsigma^{\prime}(\bu) +\alpha \nabla \cdot p \bI 
&= &\bf{f} \hspace{0.14in} \text { in } \Omega \times \mathbb{T}, \\
\bu &=&\bu_{D} \text { on } \partial \Omega_{u} \times \mathbb{T},\\
\bsigma^{\prime} \cdot \bn&=&\bsigma_{D} \text { on } \partial \Omega_{t} \times \mathbb{T}, \\
\bu&=&\bu_{0} \hspace{0.06in} \text { in } \Omega \text { at } t = 0,
\end{eqnarray}
where $\bu_D$ and ${\bsigma_D}$ are prescribed displacement and traction at boundaries, respectively, and $t$ is time.
Here, $\partial \Omega$ can be decomposed to displacement and traction boundaries, $\partial \Omega_u$ and $\partial \Omega_t$, respectively, for the solid deformation problem. 

\section{Numerical Discretizations}

For this paper, we employ the discontnious Galerkin (DG) finite element method for the spatial discretization. Let $\mathcal{T}_h$ be the shape-regular (in the sense of Ciarlet)  triangulation by a family of partitions of $\O$ into $d$-simplices $\K$ (triangles/squares in $d=2$ or tetrahedra/cubes in $d=3$). We denote by $h_{\K}$ the diameter of $\K$ and we set $h=\max_{\K \in \Th} h_{\K}$.  
Also we denote by $\Eh$ the set of all edges and by $\Eho$ and $\Ehb$ the collection of all interior and boundary edges, respectively. 
In the following notation, we assume edges for two dimension but the results hold analogously for faces in three dimensional case.
The space $H^{s}(\Th)$ $(s\in \mathbb{R})$ is the set of element-wise $H^{s}$ functions on $\mathcal{T}_h$, and $L^{2}(\Eh)$ refers to the set of functions whose traces on the elements of $\Eh$ are square integrable. Let $\mathbb{Q}_l(\K)$ denote the space of polynomials of partial degree at most $l$. 
Throughout the paper, we use the standard notation for Sobolev spaces and their norms. For example, let $E \subseteq \Omega$, then $\|\cdot\|_{1,E}$ and $|\cdot|_{1,E}$ denote the $H^1(E)$ norm and seminorm, respectively. 
For simplicity, we eliminate the subscripts on the norms if $E = \Omega$.

Since we consider the nonconforming DG methods, 
let
$$
e = \partial \K^{+}\cap \partial \K^{-}, \ \ e \in \Eho,
$$
where  $\K^{+}$ and $\K^{-}$ be two neighboring elements 
and we denote by $h_{e}$ the length of the edge $e$. 
Let $\n^{+}$ and $\n^{-}$ be the outward normal unit vectors to  $\partial T^+$ and $\partial T^-$, respectively ($\n^{\pm} :=\n_{|\K^{\pm}}$). 
For any given function $\xi$ and vector function $\bxi$, defined on the triangulation $\mathcal{T}_h$, we denote $\xi^{\pm}$ and $\bxi^{\pm}$ by the restrictions of $\xi$ and $\bxi$ to $T^\pm$, respectively. 

Next, we define the weighted average operator
$\{\cdot \}_{\delta_{e}}$  as follows:
for $\zeta \in L^2(\mathcal{T}_h)$ and $\taub \in L^2(\mathcal{T}_h)^d$,
\begin{equation}
\{\zeta\}_{\delta e}=\delta_{e} \zeta^{+}+\left(1-\delta_{e}\right) \zeta^{-}, \ \text{ and } \ 
\{\taub\}_{\delta e}=\delta_{e} \taub^{+}+\left(1-\delta_{e}\right) \taub^{-},
\ \text{ on } e \in \Eho,
\end{equation}
where $\delta_{e}$ is calculated by \cite{ErnA_StephansenA_ZuninoP-2009aa,ern2008posteriori}.
\begin{equation}
\delta_{e} :=\frac{{\kappa}^{-}_e}{{\kappa}^{+}_e+{\kappa}^{-}_e}.
\end{equation}
Here, 
\begin{equation}
{\kappa}^{+}_e :=\left(\bn^{+}\right)^{T} \cdot \bkappa^{+} \cdot \bn^{+}, \ \text{ and }
{\kappa}^{-}_e :=\left(\bn^{-}\right)^{T} \cdot \bkappa^{-} \cdot \bn^{-},
\end{equation}
where  ${\kappa_e}$ is a harmonic average of $\kappa^{+}_e$ and ${\kappa}^{-}_e$ read as:
\begin{equation}
{\kappa_{e}}:= \frac{2{\kappa}^{+}_e {\kappa}^{-}_e}{\left({\kappa}^{+}_e+{\kappa}^{-}_e\right)}.
\end{equation}
On the other hand, for $e \in \Ehb$, we set $\av{\zeta}_{\delta_e} :=   \zeta$ and $\av{\taub}_{\delta_e} :=  \taub$. 
The jump across the interior edge will be defined as 
\begin{align*}
\jump{\zeta} = \zeta^+\n^++\zeta^-\n^- \quad \mbox{ and } \quad \jtau = \taub^+\cdot\n^+ + \taub^-\cdot\n^- \quad \mbox{on } e\in \Eho. 
\end{align*}
For  $e \in \Ehb$, we let $\jump{\zeta} :=  \zeta \bn$ and $\jump{\taub} :=  \taub \cdot \bn$. 

Finally, we introduce 
the finite element space for the discontinuous Galerkin method, which is the space of piecewise discontinuous polynomials of degree $k$ by
\begin{equation}
V^{\textsf{DG}}_{h,k} (\mathcal{T}_h) := \left \{ \psi \in L^2(\Omega) | \ \psi_{|_{\K}} \in \mathbb{Q}_k(\K), \ \forall \K \in \mathcal{T}_h \right \}.
\end{equation}
Moreover, we use the notation: 
\begin{align*}
&(v,w)_{\Th}:=\dyle\sum_{\K \in \Th} \int_{\K} v\, w dx, \quad \forall\,\, v ,w \in L^{2} (\mathcal{T}_h), \\
&\langle v, w\rangle_{\Eh}:=\dyle\sum_{e\in \Eh} \int_{e} v\, w \,d\gamma, \quad \forall\, v, w \in L^{2}(\Eh).
\end{align*}

\subsection{Pressure equation}

First, we introduce the backward Euler DG approximation to \eqref{eq:mass_balance}. 
We define a partition of the time interval $0=:t^0 <t^1 < \cdots <  t^N := \mathbb{T} $ and denote the uniform time step size by $\delta t:= t^n - t^{n-1}$.
The DG finite element space approximation of the pressure $p(\bx,t)$ is denoted by $P(\bx,t) \in V^{DG}_{h,k}$. 
Let $P^n := P(\bx,t^n)$ for  $0 \leq n \leq N$. 
We set a given initial condition for the pressure as $P^0$ and assume the displacement at time $t$, $\bu(\cdot,t)$ is given.
For the simplicity the gravity and the source/sink terms are neglected.
Then, the time stepping algorithm reads as follows: Given $P^{n-1}$, 
\begin{equation}\label{eq:dgscheme}
\mbox{ Find  }P^{n} \in V_{h,k}^{DG} \mbox{  such that  } \calS_\theta(P^{n},w) = \mathcal{F}_\theta(w), \quad \forall\, w\in V_{h,k}^{DG}, \,
\end{equation}
where $\calS_\theta$ and $\mathcal{F}_\theta$ are the bilinear form and linear functional as defined by
\begin{multline}\label{stheta}  
\calS_{\theta}(v,w) :=
\dfrac{\rho}{\delta t}\left(\phi c_{f}+\frac{\alpha-\phi}{K_{s}}\right) \left ( v, w \right)_{\Th} \\
+ \left (\bkappa \nabla v,\nabla w \right )_{\Th} 
- \left \langle \av{ \bkappa \nabla v }_{\delta_{e}}, \jump{w} \right \rangle_{\mathcal{E}_h^{1}} \\
 - \theta \left \langle \jump{v},\av{ \bkappa \nabla w}_{\delta_{e}} \right \rangle_{\mathcal{E}_h^{1}} + \beta(k)  \langle h^{-1}_{e}  \kappa_e \jump{v},\jump{w} \rangle_{\mathcal{E}_h^{1}}, 
\quad \forall v, w \in V_{h,k}^{DG},  
\end{multline}
and 
\begin{multline}
\mathcal{F}_\theta(w) := \frac{1}{\delta t}(P^{n-1}, w)_{\mathcal{T}_h} 
- \rho \alpha(\dfrac{\partial}{\partial t} \nabla \cdot \bu,w)_{\mathcal{T}_h} 
- \left \langle  {q_D},\jump{w} \right \rangle_{\mathcal{E}_h^{N,\partial}} \label{ftheta} \\
\quad - \theta \left \langle p_D,\av{ \bkappa \nabla w}_{\delta_{e}} \right \rangle_{\mathcal{E}_h^{D,\partial}} 
+ \beta(k)  \langle h^{-1}_{e} \kappa_e p_D,\jump{w}\rangle_{\mathcal{E}_h^{D,\partial}},
\quad \forall w \in V_{h,k}^{DG}.
\end{multline}

 The {choice} of $\theta$ leads to different DG algorithms. For example, i) $\theta = 1$ for SIPG($\beta$)$-k$ methods \cite{StenbergR-1998aa,Dryja2003}, which later has been extended to the advection-diffusion problems in \cite{BurmanE_ZuninoP-2006aa,Di-PietroD_ErnA_GuermondJ-2008aa},
ii) $\theta = -1$ for NIPG($\beta$)$-k$ methods \cite{HoustonP_SchwabC_SuliE-2002aa}, and iii) $\theta = 0$ for IIPG($\beta$)$-k$ method \cite{DawsonC_SunS_WheelerM-2004aa}.

The interior penalty parameter, $\beta(k)$, is a function of polynomial degree approximation, $k$. Here, $h_e$ is a characteristic length of the edge $e \in \Eh$ calculated as:
\begin{equation}
h_{e} :=\frac{\operatorname{meas}\left(T^{+}\right)+\operatorname{meas}\left(T^{-}\right)}{2 \operatorname{meas}(e)},
\end{equation}
where meas($.$) represents a measurement operator, measuring length, area, or volume. Several analysis for the choice of the interior penalty parameter, $\beta$, are shown in \cite{ainsworth2007posteriori,ainsworth2010fully,ainsworth2009constant,epshteyn2007estimation,shahbazi2005explicit} and this $\beta$ is the quantatiy 
that we investigate in this paper.

\subsection{Poroelasticity problem}
\label{subsec_poro}

For the displacement $\bu$, we employ the classical continuous Galerkin (CG)finite element methods for the spatial discretizations as in \cite{choo2018enriched,Kadeethum2019} where the function space is defined as 
\begin{equation}
W^{\textsf{CG}}_{h,k} (\mathcal{T}_h) :=
\left\{\bpsi_u \in \mathbb{C}^{0}(\Omega{; \mathbb{R}^d}) :\left.\bpsi_u \right|_{T} \in \mathbb{Q}_{k}(T{; \mathbb{R}^d}), \forall T \in \mathcal{T}_{h}\right\},
\label{eq:CG_U}
\end{equation}
where
$\mathbb{C}^0(\Omega{; \mathbb{R}^d})$ denotes the space of vector-valued piecewise continuous polynomials, $\mathbb{Q}_{k}(T{; \mathbb{R}^d})$ is the space of polynomials of degree at most $k$ over each element $T$.

The CG finite element space approximation of the displacement $\bu(\bx,t)$ is denoted by $\bU(\bx,t) \in W^{CG}_{h,k}$. 
Let $\bU^n := \bU(\bx,t^n)$ for  $0 \leq n \leq N$. 
We set a given initial condition for the displacement as $\bU^0$ and the pressure at time $t$, $P^n$ is given from the previous section.
Then, the time stepping algorithm reads as follows: Given $\P^{n}$, 
\begin{equation}\label{eq:cgscheme}
\mbox{ Find  } \bU^{n} \in W_{h,k}^{CG} \mbox{  such that  } \calA(\bU^{n},\bw) = \mathcal{D}(\bw), \quad \forall\, \bw\in W_{h,k}^{CG}, \,
\end{equation}
where $\calA$ and $\mathcal{D}$ are the bilinear form and linear functional as defined as 
\begin{equation}\label{Dtheta}  
\calA(\bv,\bw) :=
\sum_{T \in \mathcal{T}_{h}} \int_{T}
\bsigma^{\prime}\left(\bv \right) : \nabla^{s} \bw \: d V  +  \sum_{T \in \mathcal{T}_{h}} \int_{T} \alpha  \nabla P^n  . \nabla \bw \: d V,
\quad \forall\, \bv, \bw\in W_{h,k}^{CG},
\end{equation}
and 
\begin{equation}
\mathcal{D}(\bw) :=
\sum_{T \in \mathcal{T}_{h}} \int_{T} {\bf{f}} \bw \: d V
+\sum_{e \in \mathcal{E}_{h}^{N}} \int_{e} \bsigma_D \bw \: d S, 
\quad \forall\, \bw\in W_{h,k}^{CG},
\end{equation}

\section{Machine Learning Algorithm}
\label{sec:alg}

\begin{figure}[!h]
\centering
\includegraphics[width=0.8\textwidth]{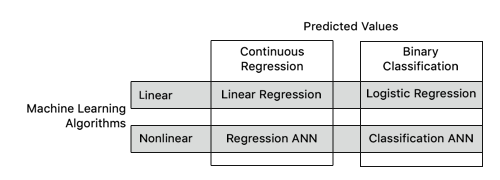}
\caption{Illustration of the methods used in this study.}
\label{fig:sec4}
\end{figure}

In this section, we present the details of the two machine learning algorithms employed in this paper to seek the effect and optimal choice of the interior penalty parameter. First, the linear approximation algorithm, which is so called linear regression or logistic regression depending on the output type, is shown. Then,  the nonlinear approximation algorithm, the artificial neural network (ANN) with deep learning algorithm, is described. See Figure \ref{fig:sec4} for the detailed outline.
Then, in the next section, the performance between linear and nonlinear approximation algorithms are compared for each given problem to find the optimal penalty parameters. 
The two different algorithms (linear/nonlinear) will provide two types of predicted values; 
one is continuous predicted value referred to as a continuous regression model, and the second one is binary predicted value referred to as a binary classification model.
Our problem to find the optimal penalty parameter could be solved by both continuous (regression model) and binary (classification model) predictions. 

First, we start by discussing the loss functions for each predicted values, continuous regression and binary classification models.
For the predicted value by employing the continuous regression model, we use mean squared error (MSE) as a loss function, which is defined as
\begin{equation}\label{eq:loss_mse}
\mathrm{MSE} :=\frac{1}{n} \sum_{i=1}^{n}\left(Y_{i}-\hat{Y}_{i}\right)^{2},
\end{equation}
where $n$ represents a number of data points, 
$Y_i$ is a true or observed values at index $i$, and $\hat{Y}_i$ is a predicted value at index $i$.  

To compare the performance of the linear regression and nonlinear ANN algorithms, by using the continuous regression predicted values,
we define $R^{2}$ and the explained variance score ($\text{EVS}$).
Here, $R^{2}$ is
\begin{equation}\label{eq:loss_r2}
R^{2} := 1-\frac{S S_{\mathrm{res}}}{S S_{\mathrm{tot}}},
\end{equation}
where $S S_{\mathrm{res}}$ is a residual sum of squares read as:
\begin{equation}\label{eq:loss_ssres}
S S_{\mathrm{res}}:=\sum_{i=1}^{n}\left(Y_{i}-\hat{Y}_{i}\right)^{2},
\end{equation}
and $S S_{\mathrm{tot}}$ is a total sum of squares defined as:
\begin{equation}\label{eq:loss_sstot}
S S_{\mathrm{tot}}:=\sum_{i=1}^{n}\left(Y_{i}-\Bar{Y}\right)^{2},
\end{equation}
where $\Bar{(\cdot)}$ is an arithmetic average operator and $\Bar{Y}$ is the arithmetic average of $Y_i$. 
Next, $\text{EVS}$ is defined as 
\begin{equation}\label{eq:loss_evs}
\text {EVS} :=1-\frac{\operatorname{Var}\{Y-\hat{Y}\}}{\operatorname{Var}\{Y\}},
\end{equation}
\noindent
where $\operatorname{Var}\{ \cdot \}$ is a variance operator. Note that if $\overline{Y-\hat{Y}} = 0$, then 
$R^{2}$ = $\text{EVS}$ or we have unbiased estimator \cite{walpole1993probability}. 

\begin{figure}[!ht]
\centering
\includegraphics[width=0.6\textwidth]{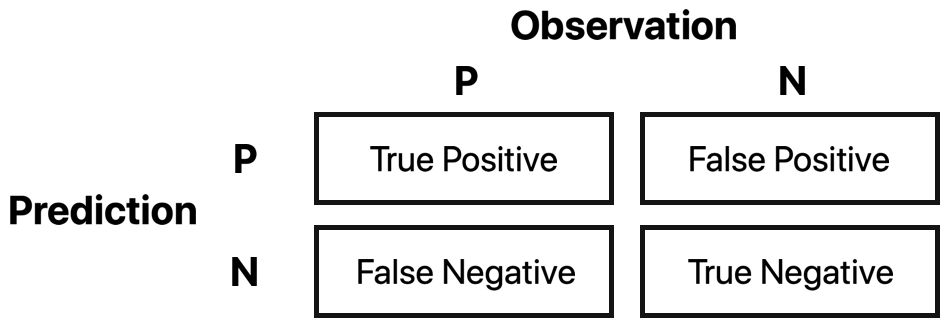}
\caption{Illustration of a confusion matrix, where \textbf{P} denotes positive and \textbf{N} denotes negative.}
\label{fig:confusion_matrix_example}
\end{figure}

For the predicted value by employing the binary classification models, 
we use binary cross entropy ($\text{BCE}$) as the loss function, which is defined as follows:
\begin{equation}\label{eq:loss_bce}
\text{BCE }:=-\frac{1}{n} 
\sum_{i=1}^{n} \left(Y_{i} \cdot \log \left(\mathbb{P}\left(Y_{i}\right)\right)+\left(1-Y_{i}\right) \cdot \log \left(1-\mathbb{P}\left(Y_{i}\right)\right)\right),
\end{equation}
\noindent
where $\mathbb{P}\left( \cdot \right)$ is a probability function. Then, we use the accuracy function ($\text{ACC}$) to compare the results from logistic regression and classification ANN, and it is defined as
\begin{equation}\label{eq:loss_acc}
\text { ACC }:=\frac{\sum \text { True positive }+\sum \text { True negative }}{n},
\end{equation}
\noindent
where `True positive' and `True negative' represent cases where the prediction agrees with the observation. See Figure \ref{fig:confusion_matrix_example} for more details.

Moreover, we employ different optimization algorithms to minimize each loss functions for linear and nonlinear algorithms, and we describe these in next sections.

\newcommand\MyBox[2]{
  \fbox{\lower0.75cm
    \vbox to 1.7cm{\vfil
      \hbox to 1.7cm{\hfil\parbox{1.4cm}{#1\\#2}\hfil}
      \vfil}%
  }%
}

\subsection{Linear approximation algorithm}
\label{sec:alg_1}

Two types of linear approximation algorithms are used in this work: (i) (multivariate) linear regression and (ii) (multivariate) logistic regression. 
These two models produce continuous and binary predictive values, respectively, 
and can consider any number of input values.
The main idea of these models is to map the linear relationship between multiple independent variables (input) 
and one dependent variable (output). 
As discussed previously, equations \eqref{eq:loss_mse} and \eqref{eq:loss_bce} are used as the loss function for multivariate linear and logistic regressions, respectively. 
To minimize these functions, we follow the classical stochastic gradient descent (SGD) \cite{scikit-learn} solver to minimize equation \eqref{eq:loss_mse}, and limited-memory Broyden-Fletcher-Goldfarb-Shanno (L-BFGS) \cite{scipy2001} solver to minimize equation \eqref{eq:loss_bce}.

Finally, we split our available data points into two parts; training set and test set using Scikit-learn package, an open-source software machine learning library for the Python programming language \cite{scikit-learn}. 
We use the same training set to train both linear and nonlinear approximation functions. 
The test set is utilized for comparing performances between the linear and nonlinear approximation algorithms. 
The splitting ratio of the data sets used in this paper is 0.8 of the total available data for the training set and 0.1 of the total available data for the test set.
We note that we only utilize 90\% of the available data set for training and test set to be consistent with the number of the data sets for the nonlinear algorithm. The nonlinear algorithm requires 10\% of the available data set for the validation set. 
Each variable input is transformed into a numeric variable, and each continuous data is normalized by its mean and variance using the prepossessing library of Scikit-learn package \cite{scikit-learn}.

\subsection{Nonlinear approximation algorithm}
\label{sec:alg_2}

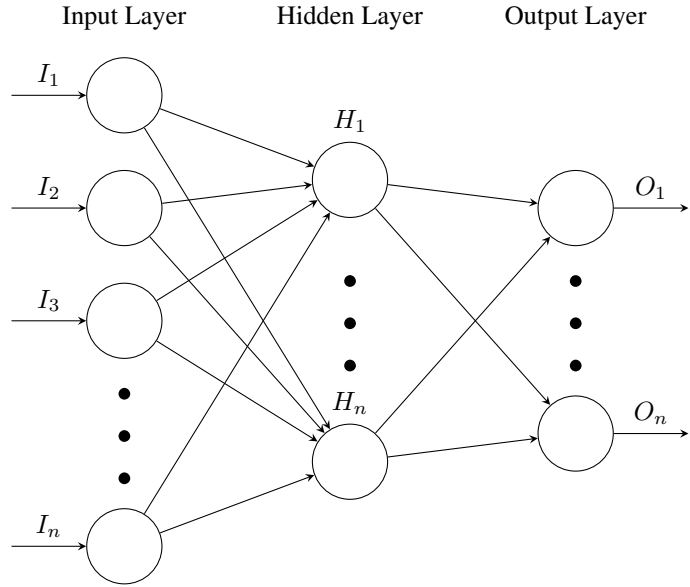
\begin{figure}[!ht]
   \centering
        \tikzset{%
  every neuron/.style={
    circle,
    draw,
    minimum size=1cm
  },
  neuron missing/.style={
    draw=none, 
    scale=4,
    text height=0.333cm,
    execute at begin node=\color{black}$\vdots$
  },
}
\begin{tikzpicture}[x=1.5cm, y=1.5cm, >=stealth]

\foreach \m/\l [count=\y] in {1,2,3,missing,4}
  \node [every neuron/.try, neuron \m/.try] (input-\m) at (0,2.5-\y) {};

\foreach \m [count=\y] in {1,missing,2}
  \node [every neuron/.try, neuron \m/.try ] (hidden-\m) at (2,2-\y*1.25) {};

\foreach \m [count=\y] in {1,missing,2}
  \node [every neuron/.try, neuron \m/.try ] (output-\m) at (4,1.5-\y) {};

\foreach \l [count=\i] in {1,2,3,n}
  \draw [<-] (input-\i) -- ++(-1,0)
    node [above, midway] {$I_\l$};

\foreach \l [count=\i] in {1,n}
  \node [above] at (hidden-\i.north) {$H_\l$};

\foreach \l [count=\i] in {1,n}
  \draw [->] (output-\i) -- ++(1,0)
    node [above, midway] {$O_\l$};

\foreach \i in {1,...,4}
  \foreach \j in {1,...,2}
    \draw [->] (input-\i) -- (hidden-\j);

\foreach \i in {1,...,2}
  \foreach \j in {1,...,2}
    \draw [->] (hidden-\i) -- (output-\j);

\node [align=center, above] at (0,2) {Input Layer};
\node [align=center, above] at (2,2) {Hidden Layer};
\node [align=center, above] at (4,2) {Output Layer};

\end{tikzpicture}
   \caption{
General neural network architecture used in this study. Input layer contains up to $i$ input nodes, and output layer is composed of $1,...,k$ output nodes. The number of hidden layers we denote $N_{hl}$ and each hidden layer is composed of $N_n$ neurons.}
   \label{fig:gerneral_nn}
\end{figure}

Similar to the previous section, we have two types of nonlinear approximation algorithms used in this work: (i) regression ANN model and (ii) classification ANN model. These models map the nonlinear relationship between input (features) to output by using nonlinear activation functions such as Sigmoid, Tanh, or rectified linear unit (ReLU) functions \cite{deng2014deep,lecun2015deep,goodfellow2016deep}. The number of hidden layers also play an important role in defining whether the neural network has deep (number of hidden layers is greater than one) or shallow (number of hidden layers is one) architecture. \cite{mhaskar2016deep}. Moreover, the shallow and deep neural networks can be combined, which is called Wide and Deep Learning, to optimize the performance and generalization \cite{cheng2016wide}. The neural network architecture used in this study is presented in Figure \ref{fig:gerneral_nn}. The number of the output node is always one, but the number of input nodes is determined from the natures of each problem, which will be discussed later. 
Hyperparameters \cite{goodfellow2016deep} are determined by 
the number of hidden layers ($N_{hl}$) and the number of neurons ($N_n$).

The artificial neural network used in this study is built on the TensorFlow platform with Keras wrapper 
\cite{tensorflow2015-whitepaper,chollet2015keras}. The ReLU is employed as the activation function for each neuron in each hidden layer for both regression and classification ANN. 
The output layer of the classification ANN is subjected to the Sigmoid activation function, while the output layer of the regression ANN is not subjected to any activation functions since the output values are continuous. 

To minimize the loss functions,  
equations \eqref{eq:loss_mse} (for the regression ANN model) and \eqref{eq:loss_bce} (for the classification ANN model), 
the mini-batch gradient descent method is used with a batch size of 10 \cite{hinton2012neural,ruder2016overview}. This method is effective since (i) it requires less memory and becomes more effective when a size of data is large and (ii) it helps to prevent gradient decent optimizations trapped in the local minimum \cite{bottou2010large,ge2015escaping}. 
Then, adaptive moment estimation (ADAM) \cite{kingma2014adam} solver is employed. 

Similar to the linear approximation algorithms, we split our available data points into three parts (i) training set, (ii) validation set, and (iii) test set using Scikit-learn package \cite{scikit-learn}.
 We note that the same training and test set data points are applied for both linear and nonlinear algorithms, but the validation set is only for the nonlinear algorithms.
The validation set is used to tune hyperparameters. 
The splitting ratio used in this paper is 0.8 of the total available data for the training set, 0.1 of the total available data for the validation set, and 0.1 of the total available data for the test set.
Some different studies suggest that it may be possible to use only training and test set and neglect the validation set when we have limited data \cite{james2013introduction}. 
In that case, hyperparameters are tuned by the test set, where the test set is also used to compare the performance. We do not prefer this case since the test set should be kept separated for the last process, and the ANN algorithms should have no prior knowledge before being tested \cite{kuhn2013applied,russell2016artificial}.

\section{Numerical results}
\label{sec:experiments}

In this section, we present several numerical experiments to illustrate the capability of our proposed algorithm and computational framework. 
All the numerical experiments are computed by the program built by employing Scikit \cite{scikit-learn} for linear algorithms and TensorFlow with Keras wrapper \cite{tensorflow2015-whitepaper,chollet2015keras} for nonlinear algorithms. Moreover, 
the presented results are also computed by the JMP platform (SAS) \cite{thomas1999jmp} to verify our results. 
Mainly, the continuous regression for the predicted values is used for Section \ref{subsec_elliptic} to find the optimal penalty parameter for the elliptic problem in each case. 
Thus, the performance of the linear regression and nonlinear regression ANN is compared in this section. 
For Section \ref{subsec_biot}, to find the optimal penalty parameter for the Biot's equation in each case, the binary classification for the predicted values is used. Here, the performance of the linear logistic regression and nonlinear classification ANN methods is compared.

\subsection{Effect and optimal choice of interior penalty parameter for elliptic equations}
\label{subsec_elliptic}

First, we study the effect  and optimal choice of the interior penalty parameter $\beta$ 
by considering the simplified elliptic equation of the flow problem  \eqref{eq:mass_balance}.
By assuming that the pressure does not depend on the time and $\bu$ is a given constant, we obtain the simplified equation
\begin{equation} \label{eq:mass_balance_elliptic}
-\nabla \cdot ( \bkappa \nabla p) = g \text { in } \Omega. 
\end{equation}
We note that the gravity ($\bg$) is  neglected for simplicity, and $g$ is the source/sink term. 
In each of the following problems, we compare the performance of the linear regression and nonlinear regression ANN, where the predicted values are continuous.

\subsubsection{The effect of a polynomial degree approximation ($k$)} 
\label{subsubsec_1}
Before we employ the machine learning algorithm, presented in Section \ref{sec:alg}, we investigate the effect of the polynomial degree approximation to the penalty parameter. We illustrate the effect of a polynomial degree approximation on the choice of optimal $\beta$ using different linear solvers (direct or iterative solver) and discretization schemes. Here different discretization schemes indicate the options for choosing IIPG ($\theta = 0$) or SIPG ($\theta = 1$).

For this case, we set the exact solution in $\Omega = \left[ 0,1\right]^2$ as
\begin{equation}\label{eq:darcy_exact_solution}
p(x,y) := \sin(x+y),
\end{equation}
and $\kappa$ ($\boldsymbol{\kappa} := \kappa  \boldsymbol{I}$) has different values in a range of $[1.0 \times 10^{-18},1.0]$. 
Furthermore, the homogeneous boundary conditions are applied to all boundaries.
In particular, we study the five different $k$ values (1, 2, 3, 4, and 5), and each cases are tested by different combination of linear solvers and $\theta$. The detailed algorithm is presented in Algorithm \ref{alg:ell_effect_of_k}.

\begin{algorithm}[!ht]
\caption{Investigation procedure for the elliptic problem}
\label{alg:ell_effect_of_k}
\begin{algorithmic}

\STATE{Initialize the data set of ${\kappa}$ that used in the investigation}
\FOR{$i$ $<$ $n_{{\kappa}}$, where $n_{\boldsymbol{\kappa}}$ is the size of the specified data set ${\kappa}$,}
\STATE{Assign $\boldsymbol{\kappa}:=\kappa\left[i\right]\boldsymbol{I}$, $\beta:=\beta_0$, 
where $\beta_0 = 100.00$,
and $h:=h_0$, where $h_0 = 6.25 \times 10^{-2}$}
\WHILE{error convergence rate is optimal}
\STATE{Update $\beta:=0.99 \times \beta$} \COMMENT{Except the first loop}
\FOR{$j$ $<$ $n_h$, where $n_h = 6$}
\STATE{Solve \eqref{eq:mass_balance_elliptic}
and compute the error}
\STATE{Calculate error convergence rate}
\STATE{Update $h:=0.5 \times h$}
\ENDFOR
\ENDWHILE
\RETURN $\beta$, this $\beta$ is the smallest, which the optimal error convergence rate can be observed.
\ENDFOR
\end{algorithmic}
\end{algorithm}

The main idea of this algorithm is that we reduce the $\beta$ values (1\% by each test) 
until the optimal error convergence rate is not guaranteed anymore. Thus, in other words, we focus on finding the smallest $\beta$ that ensures the optimal convergence rate.
Here, the optimal convergence rate is obtained by six computations cycle on uniform triangular meshes, where the mesh size $h$ is divided by two for each cycle. The behavior of the $H^1(\Omega)$ semi norm errors for the approximated solution versus the mesh size $h$ is checked.

The results presented in Table \ref{tab:ell_homo_beta_unstable} show that the lowest (optimal) $\beta$ values for SIPG, 
($\theta = 1$) for each cases, are higher than those for IIPG ($\theta = 0$). 
Besides, the smallest $\beta$ values increase as $k$ increases. 
However, the choice of linear solver, either direct or iterative solver, did not influence the results. 
Theses computations are implemented by using  FEniCS \cite{AlnaesBlechta2015a}, and the direct solver used in this problem is lower-upper decomposition (LU) while conjugate gradient (CG) method with algebraic multigrid methods (AMG) method preconditioner \cite{petsc-user-ref} are employed for an iterative scheme. 
\begin{table}[!h]
  \centering
    \begin{tabular}{|c|c|c|c|c|}
    \hline
    \multirow{2}[4]{*}{$k$} & \multicolumn{2}{c|}{SIPG ($\theta=1$)} & \multicolumn{2}{c|}{IIPG  ($\theta=0$)} \\
\cline{2-5}          & \multicolumn{1}{l|}{direct solver} & \multicolumn{1}{l|}{iterative solver} & \multicolumn{1}{l|}{direct solver} & \multicolumn{1}{l|}{iterative solver} \\
    \hline
    1     & 1.11  & 1.11  & 0.83  & 0.83 \\
    \hline
    2     & 2.80  & 2.80  & 2.74  & 2.74 \\
    \hline
    3     & 5.79  & 5.79  & 5.68  & 5.68 \\
    \hline
    4     & 9.99  & 9.99  & 9.79  & 9.79 \\
    \hline
    5     & 14.97 & 14.97 & 14.67 & 14.67 \\
    \hline
    \end{tabular}%
    \caption{The lowest $\beta$ value that provides the optimal error convergence rate solution with different $k$, $\theta$, and linear solver. 
    }
  \label{tab:ell_homo_beta_unstable}%
\end{table}%

\subsubsection{Effect of interior penalty parameter for linear solvers and optimal choice by employing machine learning algorithms}
\label{subsubsec_512}

\begin{figure}[!h]
   \centering
        \includegraphics[width=0.25\textwidth]{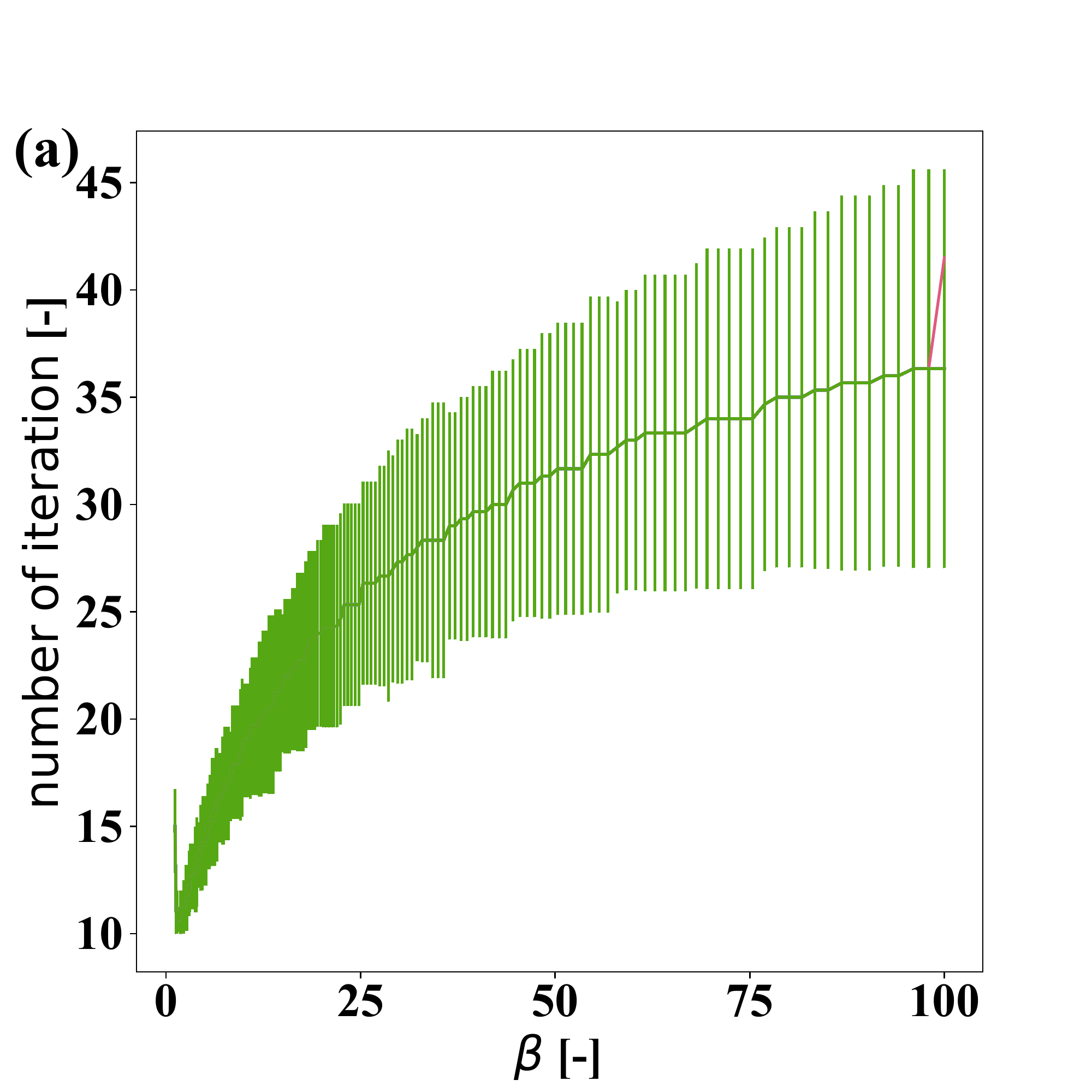}
        \includegraphics[width=0.25\textwidth]{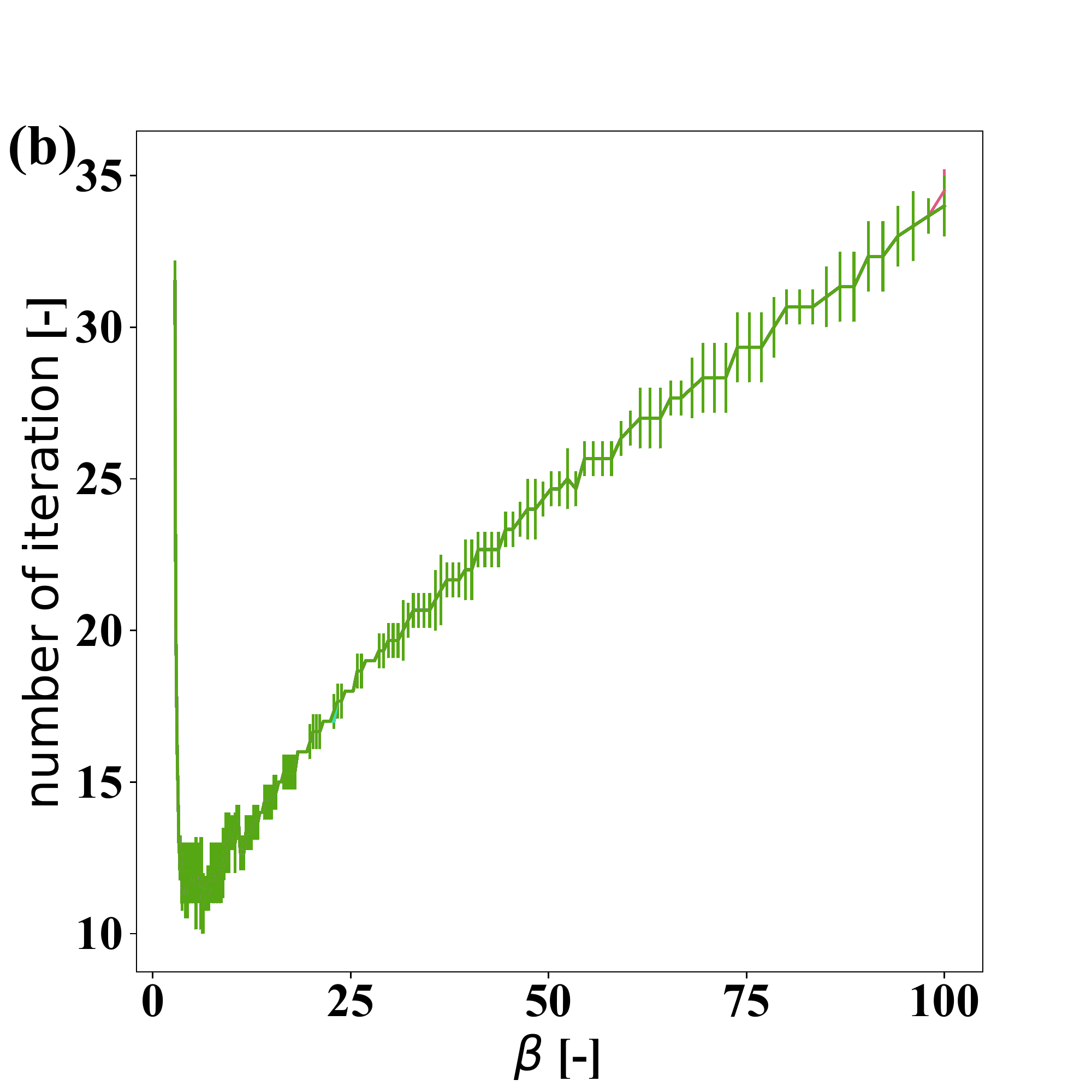}
        \includegraphics[width=0.25\textwidth]{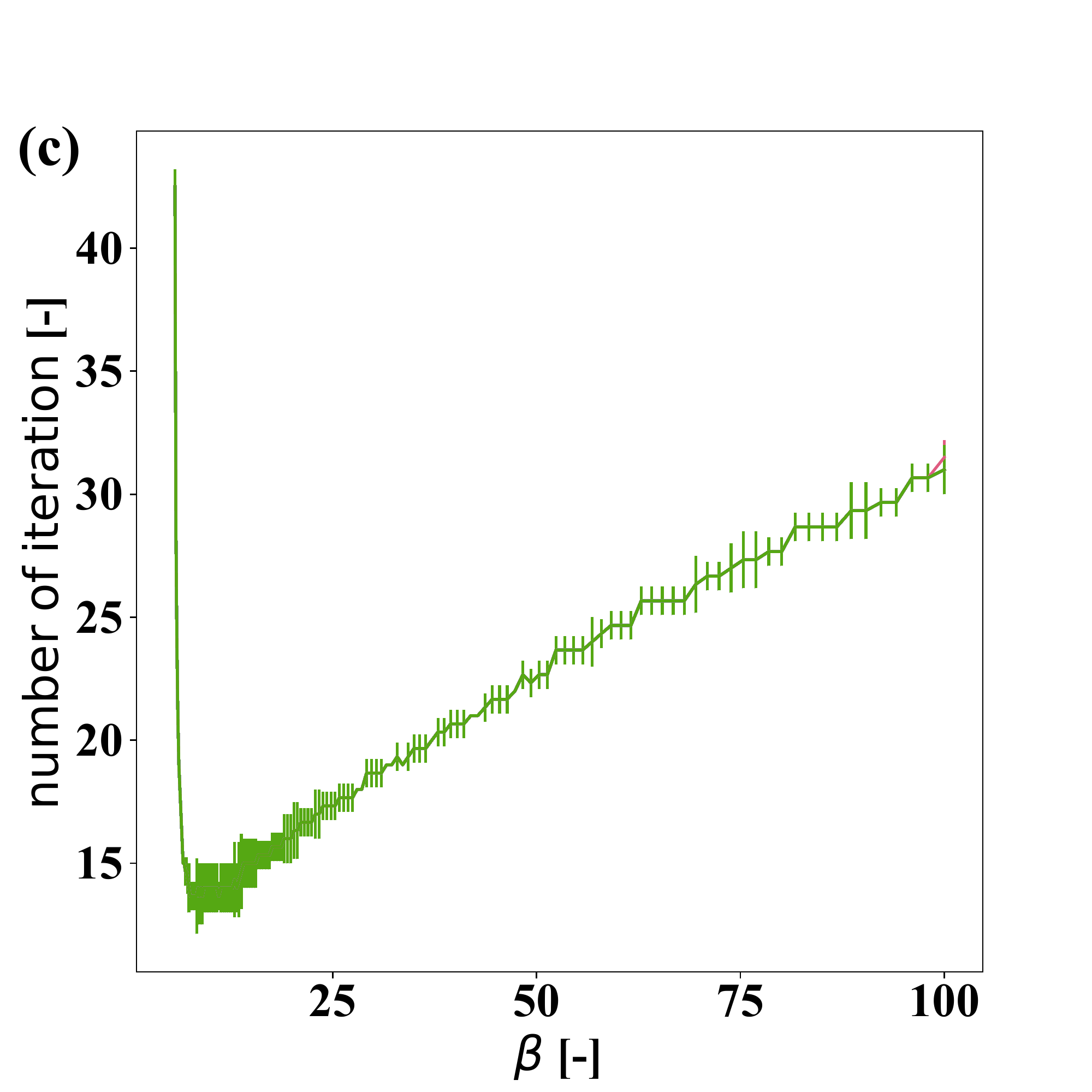}
        \includegraphics[width=0.25\textwidth]{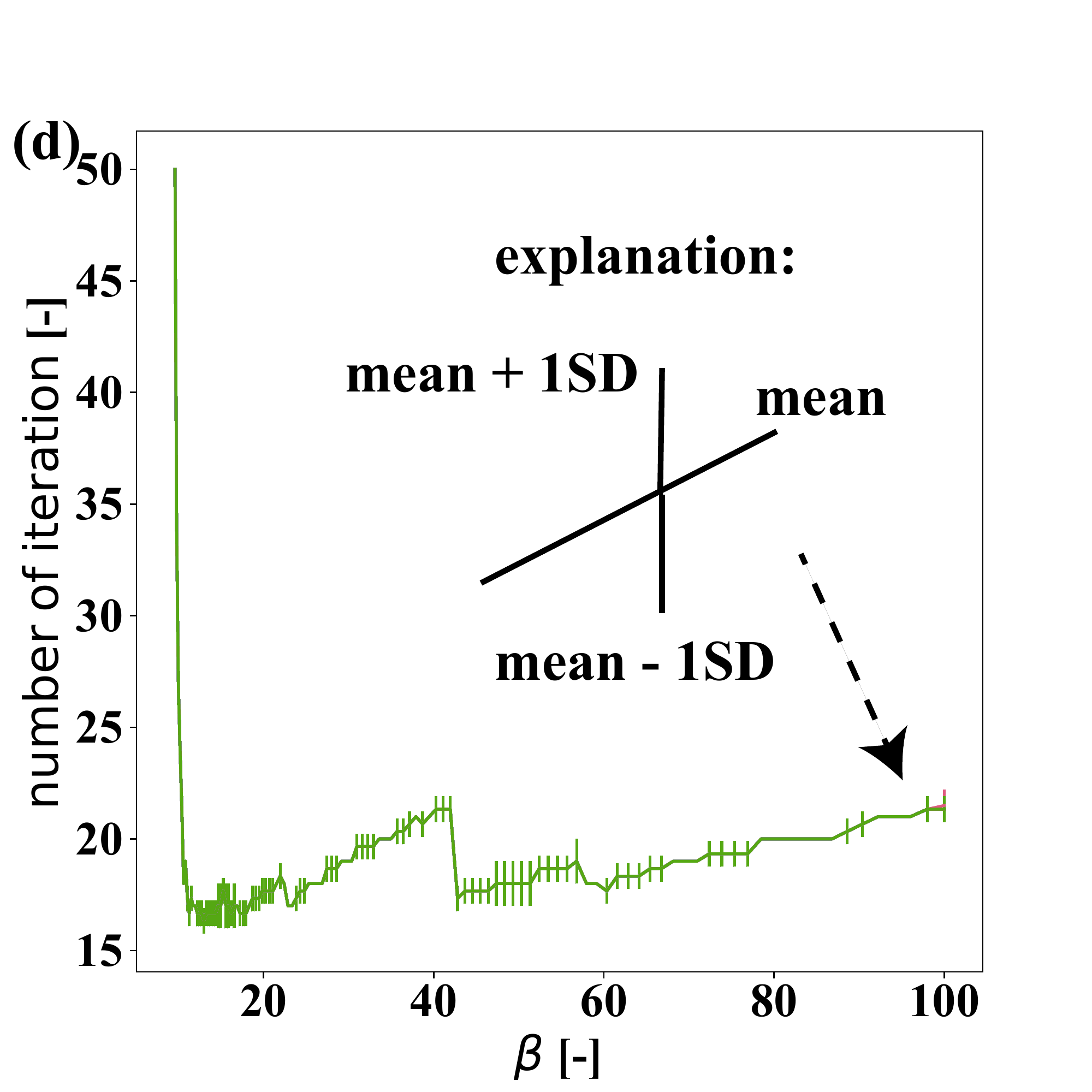}
        \includegraphics[width=0.25\textwidth]{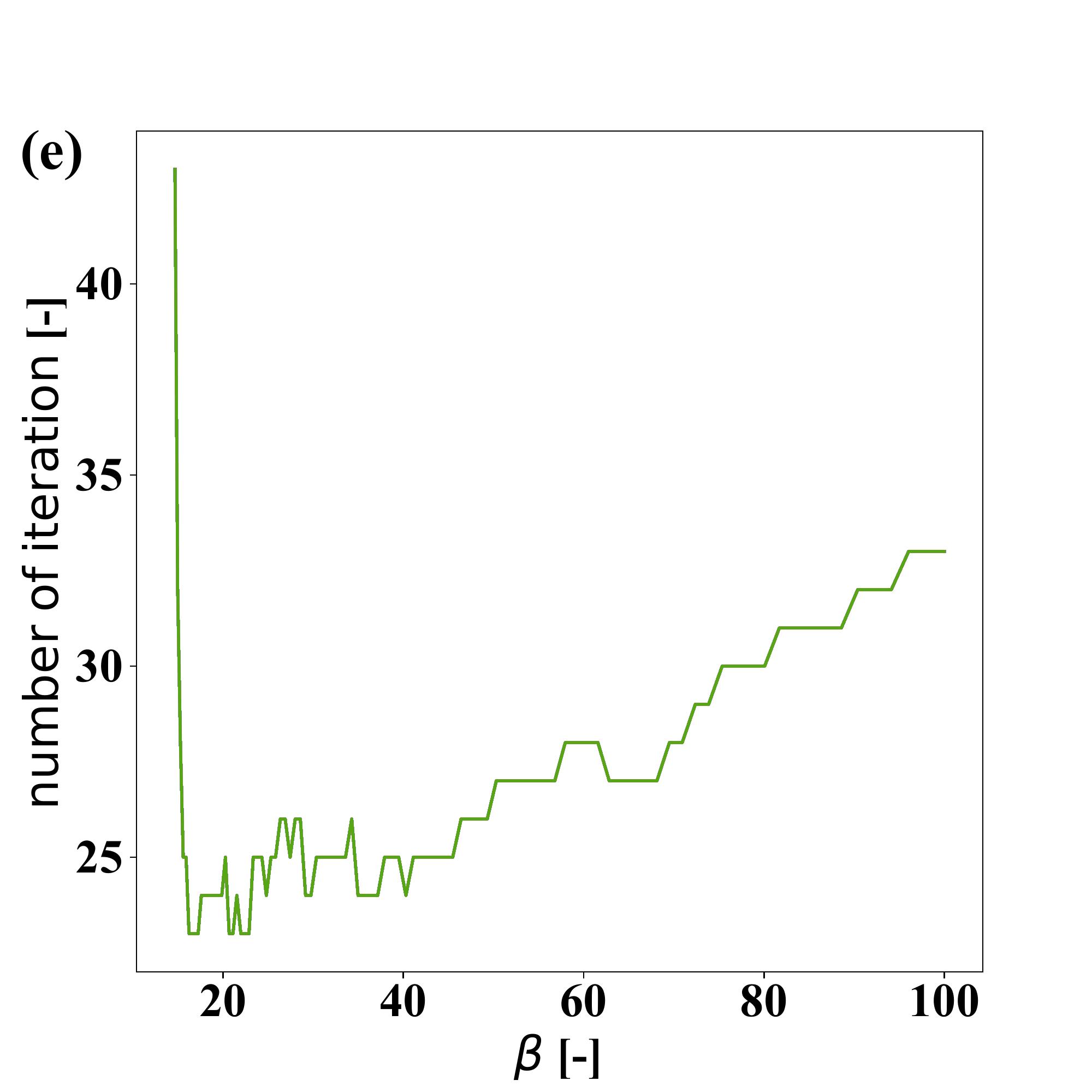}
\caption{Number of linear iterative solver of SIPG for (\textbf{a}) $k=1$, (\textbf{b}) $k=2$, (\textbf{c}) $k=3$, (\textbf{d}) $k=4$, and (\textbf{e}) $k=5$. Note that each line represents different value of ${\kappa}$, and the error bar shows mean and standard derivation ($\pm$ 1 SD) of each number of iteration see (\textbf{d}) for an explanation.}
   \label{fig:ell_homo_sipg_noi}
\end{figure}

\begin{figure}[!h]
   \centering
        \includegraphics[width=0.25\textwidth]{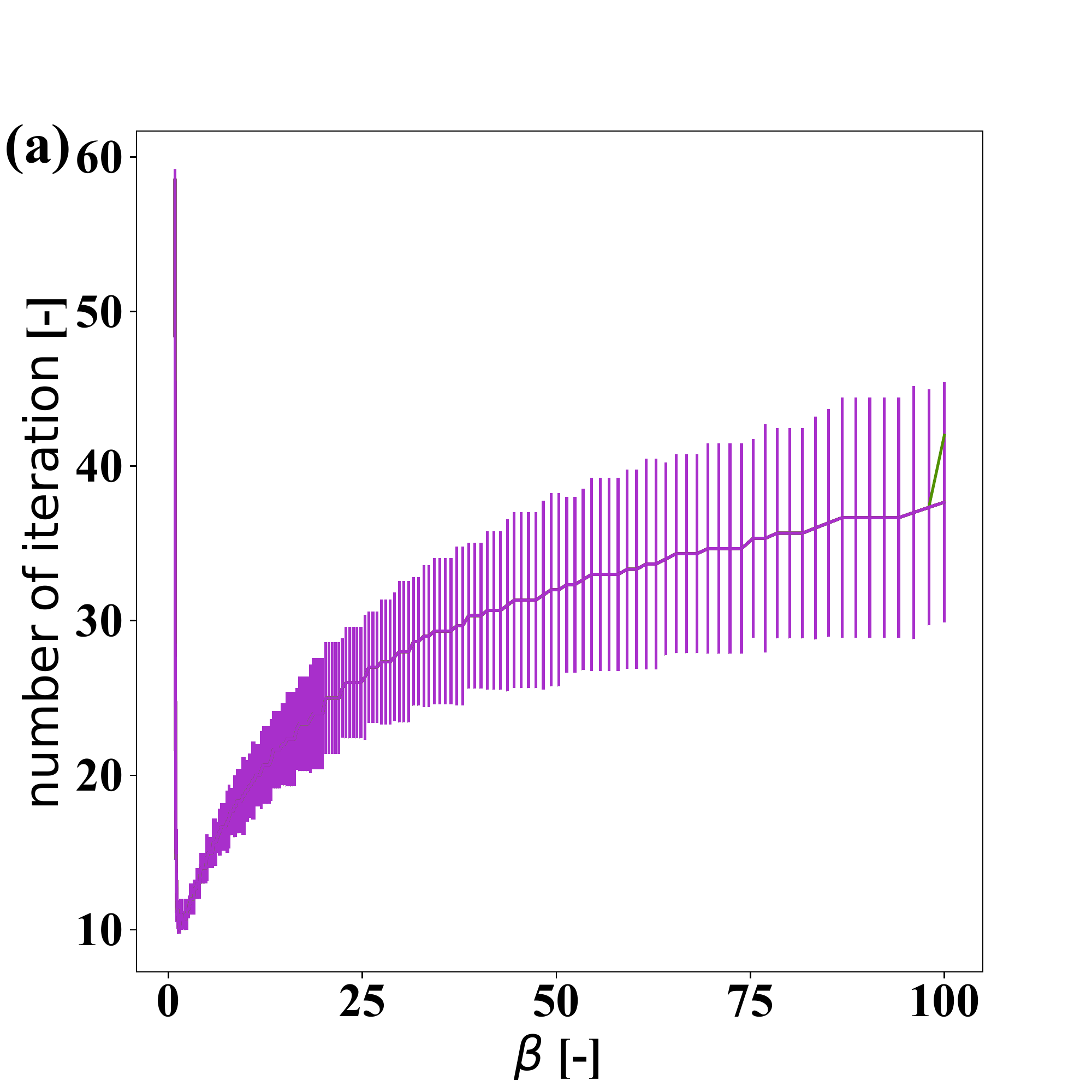}
        \includegraphics[width=0.25\textwidth]{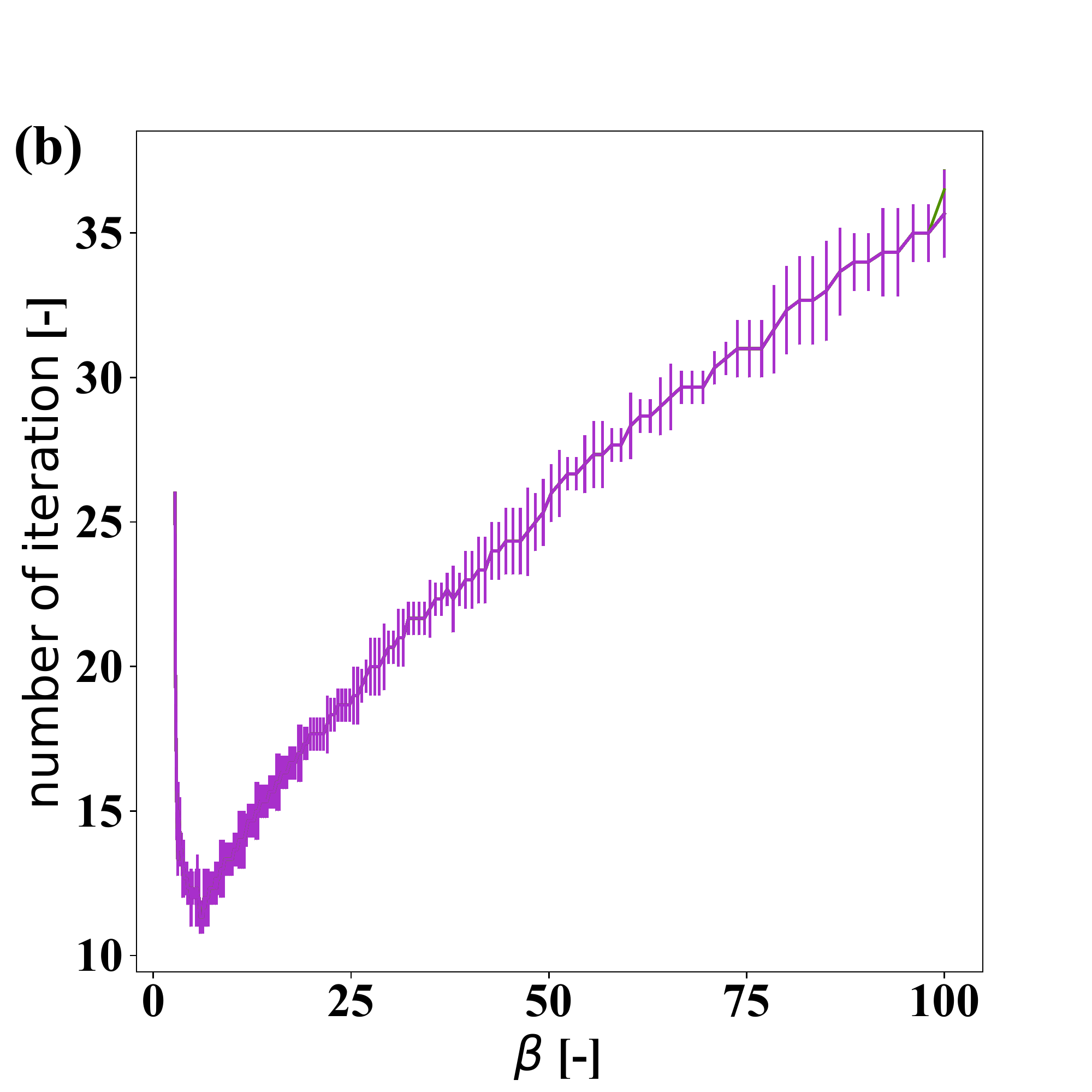}
        \includegraphics[width=0.25\textwidth]{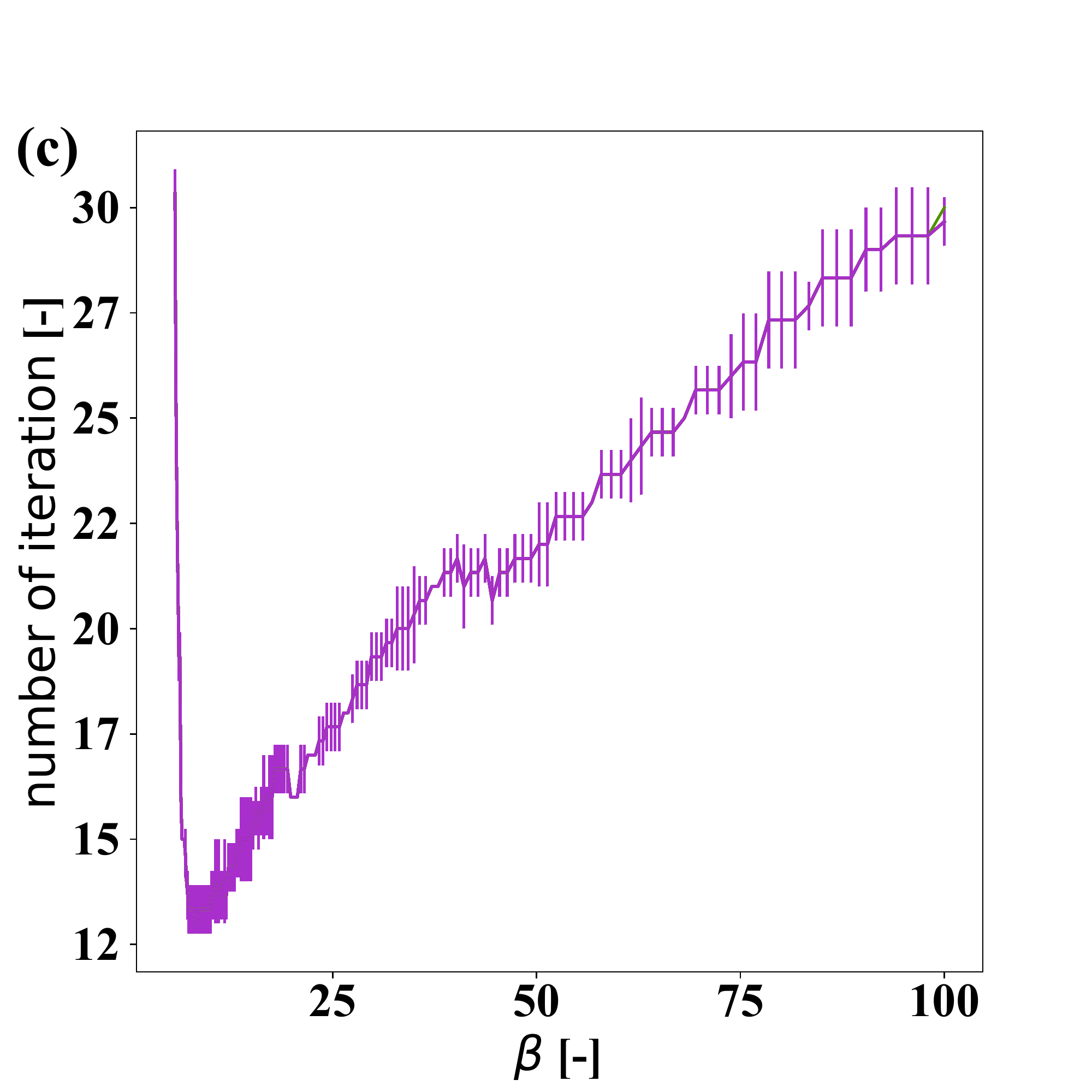}
        \includegraphics[width=0.25\textwidth]{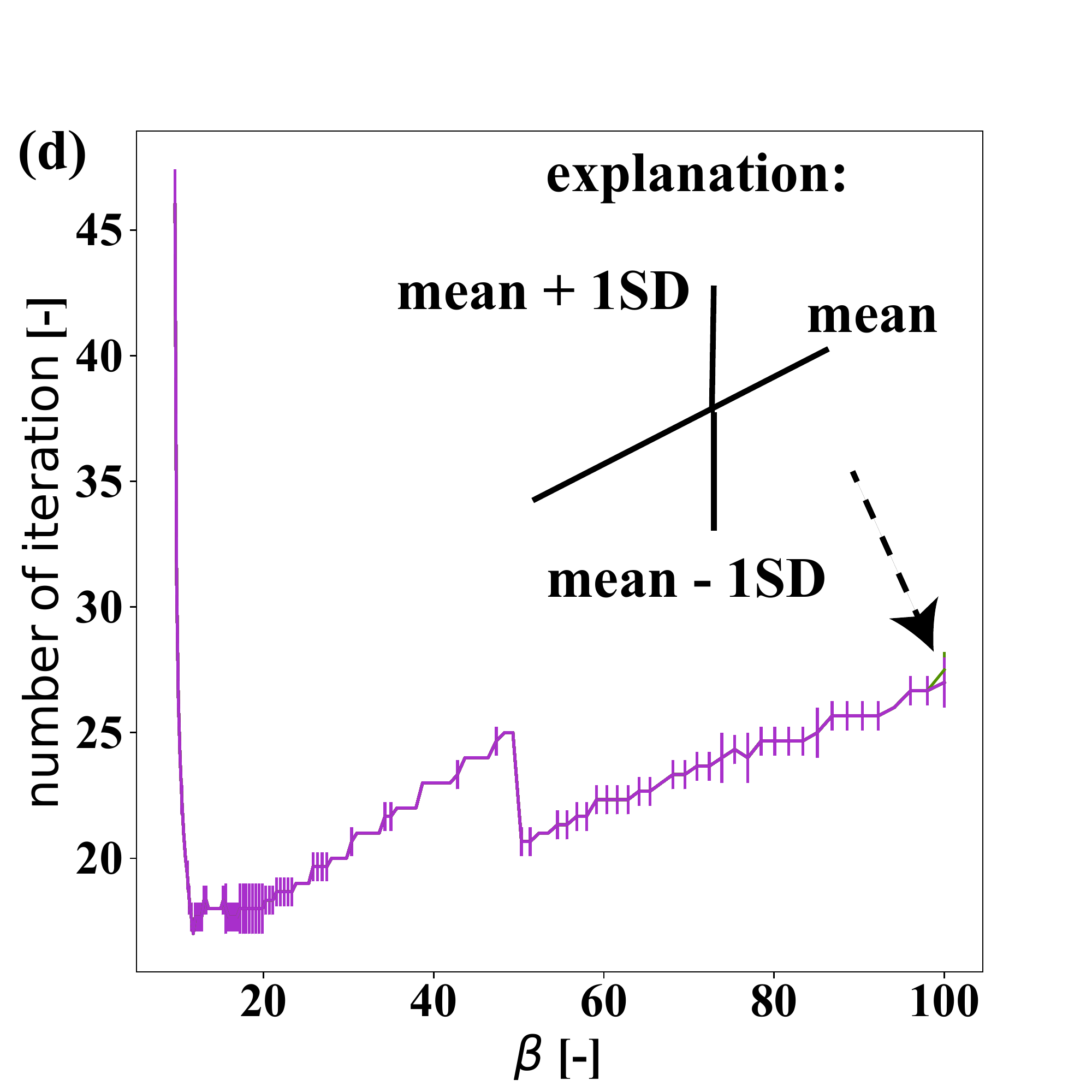}
        \includegraphics[width=0.25\textwidth]{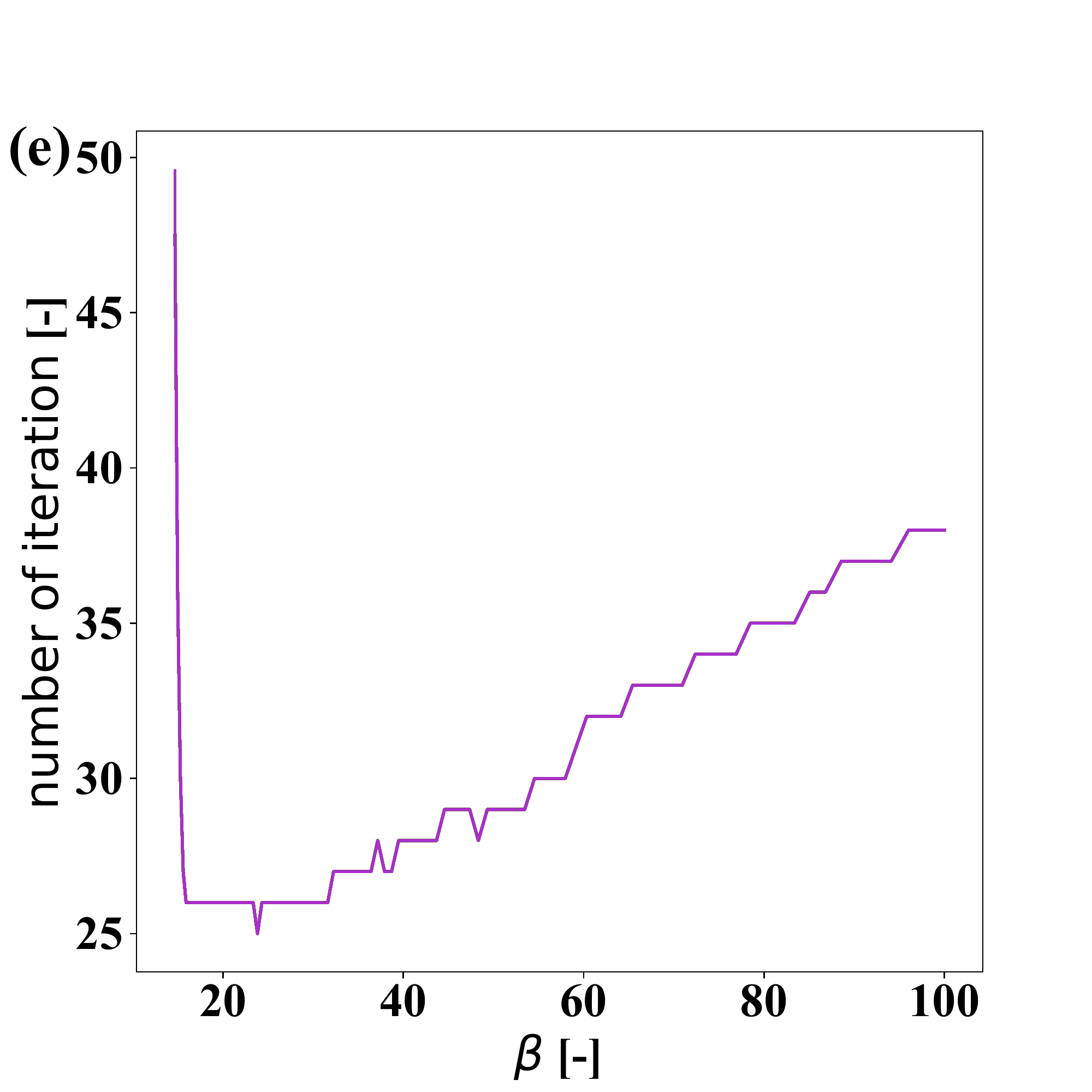}
   \caption{Number of linear iterative solver of IIPG for (\textbf{a}) $k=1$, (\textbf{b}) $k=2$, (\textbf{c}) $k=3$, (\textbf{d}) $k=4$, and (\textbf{e}) $k=5$. Note that each line represents different value of ${\kappa}$, and the error bar shows mean and standard derivation ($\pm$ 1 SD) of each number of iteration see (\textbf{d}) for an explanation.}
   \label{fig:ell_homo_iipg_noi}
\end{figure}

However, it is observed that the choice of $\beta$ values significantly impacts the number of iterations for the linear solver.  
Figures \ref{fig:ell_homo_sipg_noi} and \ref{fig:ell_homo_iipg_noi} illustrate the number of iterations of the linear solver for SIPG and IIPG, respectively.
In the beginning, the number of iteration decreases, when $\beta$ decreases, but
when $\beta$ approaches zero, the number of iteration increases dramatically. 
Subsequently, the solver becomes unstable and doesn't converge to the solution.


Thus, we confirmed that the choice of $\beta$ is essential for the linear solvers. 
To be precise, if $\beta$ is too large, the iteration number is high, but also the too small value of $\beta$ can cause a high number of iteration and, more importantly, non-convergence.

To find the optimal $\beta$ for the iterative solver, we need to consider $\beta$, which requires a minimum of the linear solver iteration, provides stable solutions, and optimal error convergence rate is ensured. 
Hence, we first identify the parameters, which impact the number of iteration (dependent variable) by employing the chi-squared test \cite{scikit-learn}.  Table \ref{tab:linear_regression_p_value} illustrates the results of the test and p-values for each variables are presented.  
We note that $\theta$, $\beta$, $h$, and $k$ have a p-value of less than 0.025; therefore, we include these variables as independent variables for further predictive model development. 
The $\kappa$, however, does not affect the results since $\kappa$ is included in  a coefficient of the penalty term as shown in \eqref{stheta}.

\begin{table}[!h]
  \centering
    \begin{tabular}{|c|c|}
    \hline
    Variable & p-value \\
    \hline
    $\theta$ & $\approx$ 0.00 \\
    \hline
    $\kappa$ & $\approx$ 1.00 \\
    \hline
    $\beta$  & $\approx$ 0.00 \\
    \hline
    $h$     & $\approx$ 0.00 \\
    \hline
    $k$     & $\approx$ 0.00 \\
    \hline
    \end{tabular}%
\caption{p-value results for each explanatory variables for elliptic equation.}
\label{tab:linear_regression_p_value}%
\end{table}%

Subsequently, from the results in Table \ref{tab:linear_regression_p_value}, we employ the linear and nonlinear machine learning algorithms that were presented in Section \ref{sec:alg} to find the optimal choice of $\beta$ which ensures both the minimum iteration number for the linear solver and optimal convergence rate (stability). To elaborate, we want to find a range of $\beta$ that could guarantee the stability, see Table \ref{tab:ell_homo_beta_unstable}, while utilizing the minimum of number of iteration, see Figures \ref{fig:ell_homo_sipg_noi} and \ref{fig:ell_homo_iipg_noi}.

In this problem, we have total 182,881 data sets (all the values we plot on Figures \ref{fig:ell_homo_sipg_noi} and \ref{fig:ell_homo_iipg_noi}). 
As discussed in both sections \ref{sec:alg_1} and \ref{sec:alg_2}, the data sets are split by  training, validation, and test sets  
using the splitting ratio 
$\left[ 0.8, 0.1 ,0.1\right]$. Thus, the number of  training sets, validation sets, and test sets are  $0.8 \times 182,881$, $0.1 \times 182881$, and  $0.1 \times 182881$, respectively.
We use the training set to train the linear and nonlinear machine learning algorithms. The validation set is for tuning the hyperparameters for the nonlinear ANN models, and the test set is for comparing performances between the linear and nonlinear algorithms.

First, we begin with the linear algorithm by building the multi-variable regression \cite{scikit-learn} as follows:
\begin{equation} \label{eq:ell_regression_con}
\begin{split}
\mathrm{number \ of \ iteration} = \alpha 
    + \gamma_{1} \times \theta 
    + \gamma_{2} \times \beta    
    + \gamma_{3} \times h 
    + \gamma_{4} \times k
\end{split}
\end{equation}
where $\alpha = 16.93$, $\gamma_{1} = -1.11$,  $\gamma_{2} = 0.21$, $\gamma_{3} = -9.85$, and $\gamma_{4} = 0.02$. These parameters provide the minimum value ($\leq 1 \times 10^{-4}$) of MSE value \eqref{eq:loss_mse}.  
Then, we obtain the $r^2$ and explained variance score ($\text{EVS}$) as 
\begin{equation}
r^2 = 0.60 \ \text{ and } \ \text{ EVS } = 0.60,
\label{results_linear_1}
\end{equation}
\noindent
by \eqref{eq:loss_r2} and \eqref{eq:loss_evs} as expalined in section \ref{sec:alg}.

Secondly, to compare the above linear algorithm with the nonlinear ANN algorithm, we construct the nonlinear ANN algorithm by using four inputs ($\theta,\beta,h$, and $k$ and one output (number of iteration) as presented in Figure \ref{fig:ell_ann_fig}.  
For simplicity, we assume each hidden layer has the same number of neuron and Rectified Linear Unit (ReLU) is used as an activation function for each neuron of the hidden layer. ADAM \cite{kingma2014adam} is used to minimize the loss function, which is MSE \eqref{eq:loss_mse} in this case. 

\tikzset{%
  every neuron/.style={
    circle,
    draw,
    minimum size=0.1cm
  },
  neuron missing/.style={
    draw=none, 
    scale=1.5,
    text height=0.333cm,
    execute at begin node=\color{black}$\vdots$
  },
}
\begin{figure}[!h]
    \centering
    \begin{tikzpicture}[x=1.5cm, y=1.5cm, >=stealth, scale=0.5]
\foreach \m/\l [count=\y] in {1,2,3,4}
  \node [every neuron/.try, neuron \m/.try] (input-\m) at (0,2.1-\y) {};

\foreach \m [count=\y] in {1,missing,2}
  \node [every neuron/.try, neuron \m/.try ] (hidden-\m) at (2,2-\y*1.25) {};

\foreach \m [count=\y] in {1}
  \node [every neuron/.try, neuron \m/.try ] (output-\m) at (4,0.75-\y) {};

\draw [<-] (input-1) -- ++(-1,0)
    node [above, midway] {$\theta$};
\draw [<-] (input-2) -- ++(-1,0)
    node [above, midway] {$\beta$};
\draw [<-] (input-3) -- ++(-1,0)
    node [above, midway] {$h$};
\draw [<-] (input-4) -- ++(-1,0)
    node [above, midway] {$k$};


\foreach \l [count=\i] in {1,n}
  \node [above] at (hidden-\i.north) {$H_\l$};

\draw [->] (output-1) -- ++(4.3,0)
    node [above, midway] {number of iteration};

\foreach \i in {1,...,4}
  \foreach \j in {1,...,2}
    \draw [->] (input-\i) -- (hidden-\j);

\foreach \i in {1,2}
  \foreach \j in {1}
    \draw [->] (hidden-\i) -- (output-\j);


\node [align=center, above] at (0,2) {\footnotesize Input Layer };
\node [align=center, above] at (2.4,2) {\footnotesize Hidden Layer };
\node [align=center, above] at (4.8,2) {\footnotesize  Output Layer};

\end{tikzpicture}
    \caption{Neural network architecture used for the elliptic problem with continuous exact solution. The number of hidden layers, $N_{hl}$, and the number of neuron for each hidden layer, $N_n$, are used as the sensitivity analysis parameters. $H_1$ and $H_n$ represent the numbering of each neuron in each hidden layer.}
    \label{fig:ell_ann_fig}
\end{figure}
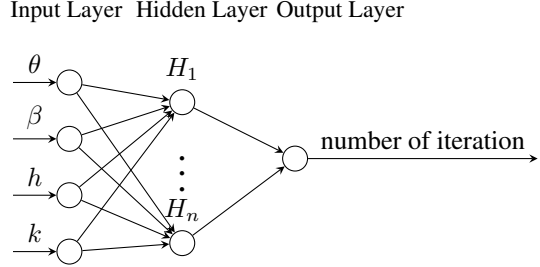

\begin{table}[!ht]
\centering
\begin{tabular}{|c|c|c|c|c|}
\hline
\backslashbox{$N_{hl}$}{$N_n$}   & 10      & 20      & 40     & 80       \\ \hline

2     & 4.27  & 1.41  & 1.20  & 0.95 \\
    \hline
    4     & 2.21  & 1.43  & 2.13  & 1.58 \\
    \hline
    8     & 3.60  & 1.74  & 1.60  & 0.98 \\
    \hline

\end{tabular}
\caption{Elliptic equation with continuous exact solution: Mean squared error (MSE) values of the validation set for different number of hidden layers $N_{hl}$ and different number of neurons per layer $N_{n}$}
\label{table:ell_hyper_nl_nn_con}
\end{table}

Table \ref{table:ell_hyper_nl_nn_con} illustrates that the MSE of the validation set is generally decreased as $N_{hl}$ and $N_{n}$ are increased.
Since we observe that the neural network performance is not significantly improved when $N_{hl} > 2$ and $N_{n} > 80$, which shows the sign of overfitting,
we choose $N_{hl} = 2$ and $N_{n} = 80$ for the test set. 
Then the final results for the nonlinear ANN algorithm give 

\begin{equation}
r^2 = 0.98 \ \text{ and } \ \text{EVS} = 0.98.
\label{results_nonlinear_1}
\end{equation}
By comparing 
\eqref{results_linear_1} 
and \eqref{results_nonlinear_1},
we note that the results from the nonlinear ANN \eqref{results_nonlinear_1} illustrate the significant improvement of the prediction performance as the $r^2$ and EVS are improved from the linear algorithm (multi-variable linear regression) \eqref{results_linear_1} significantly.

\subsubsection{The effect of the continuity of the solutions and a heterogeneous coefficient} 
Next, we investigate the effect on the choice of optimal $\beta$ by the continuity of the solutions, heterogeneity of $\boldsymbol{\kappa}$,
and $\theta$ values. 

For the continuous solution, we take the same exact solution \eqref{eq:darcy_exact_solution} as used in section \ref{subsubsec_1}. 
However, in this example, we choose the heterogeneous coefficient by setting:
\begin{equation}\label{eq:darcy_exact_solution_1_kappa}
\boldsymbol{\kappa} := \kappa  \sin(x+y) \boldsymbol{I}.
\end{equation}

Next, for the discontinuous solution, we 
set the exact solution in $\Omega  = \left[ 0,1\right]^1$ as: 
\begin{equation}\label{eq:dis_exact_solution}
p =
\begin{cases}
2.0x \dfrac{{\kappa}_1}{{\kappa}_0+{\kappa}_1} \quad & \text { if } x \leq 0.5, \vspace{0.1in} \\
\dfrac{\left( 2x-1 \right)  {\kappa}_0+{\kappa}_1}{{\kappa}_0+{\kappa}_1} \quad & \text { if } x > 0.5,
\end{cases}
\end{equation}

\noindent
where ${\kappa}_1$ and ${\kappa}_2$ represent multiplied coefficients for the $0 \leq x \leq 0.5$ and $1.0 \geq x > 0.5$ subdomains, respectively. 
Here, ${\kappa}_1$ $\neq$ ${\kappa}_2$, and ${\kappa}_1$ and ${\kappa}_2$ have the same range of $[1.0, 1.0 \times 10^{-18}]$. Then $\boldsymbol{\kappa}$ for the discontinuous solution is
\begin{equation}\label{eq:darcy_exact_solution_dis_kappa}
\boldsymbol{\kappa} :=
\begin{cases}
{\kappa}_1 \boldsymbol{I} \quad & \text { if }  0.0 \leq x \leq 0.5, \\
{\kappa}_2 \boldsymbol{I} \quad & \text { if }  1.0 \geq x > 0.5.
\end{cases}
\end{equation}
Subsequently, the boundary conditions are applied as follows:
\begin{equation}\label{eq:dis_bound_solution}
p =
\begin{cases}
0.0 \quad & \text { at } x = 0.0, \\
1.0 \quad & \text { at } x = 1.0,
\end{cases}
\end{equation}

To find the optimal $\beta$ which provides the optimal error convergence rate, we employ the Algorithm \ref{alg:ell_effect_of_k}. 
In this case, we set $k=1$ but vary the choice of linear solver, $\theta$, and the exact solutions (continuous/discontinuous). The results for the optimal $\beta$ are shown in Table \ref{tab:ell_homo_beta_unstable_dis}. The results of SIPG illustrate the similarity between the continuous and discontinuous solutions. The results of IIPG, however, show a discrepancy as to the lowest $\beta$ values that provide the optimal convergence rate solution are different between the continuous and discontinuous solutions. The type of solver, direct and iterative solvers, does not influence the results.

\begin{table}[h!]
  \centering
    \begin{tabular}{|c|c|c|c|c|}
    \hline
    \multirow{2}[4]{*}{exact solution} & \multicolumn{2}{c|}{SIPG} & \multicolumn{2}{c|}{IIPG} \\
\cline{2-5}          & \multicolumn{1}{l|}{direct solver} & \multicolumn{1}{l|}{iterative solver} & \multicolumn{1}{l|}{direct solver} & \multicolumn{1}{l|}{iterative solver} \\
    \hline
    continuous \eqref{eq:darcy_exact_solution} & \multicolumn{1}{c|}{1.11} & \multicolumn{1}{c|}{1.11} & \multicolumn{1}{c|}{0.83} & \multicolumn{1}{c|}{0.83} \\
    \hline
    discontinuous \eqref{eq:dis_exact_solution}& 1.11       &    1.11   &    0.89   & 0.89  \\
    \hline
    \end{tabular}%
    \caption{The lowest $\beta$ value that provides the optimal convergence rate solution with different type of exact solution (continuous or discontinuous), $\theta$, and linear solver. Note that $\boldsymbol{{\kappa}}$ is heterogeneous, and $k=1$.}
  
  \label{tab:ell_homo_beta_unstable_dis}%
\end{table}%

\subsubsection{Effect of interior penalty parameter for linear solvers and optimal choice by employing machine learning algorithms}
\label{subsubsec_514}

\begin{figure}[h!]
   \centering
        \includegraphics[width=0.3\textwidth]{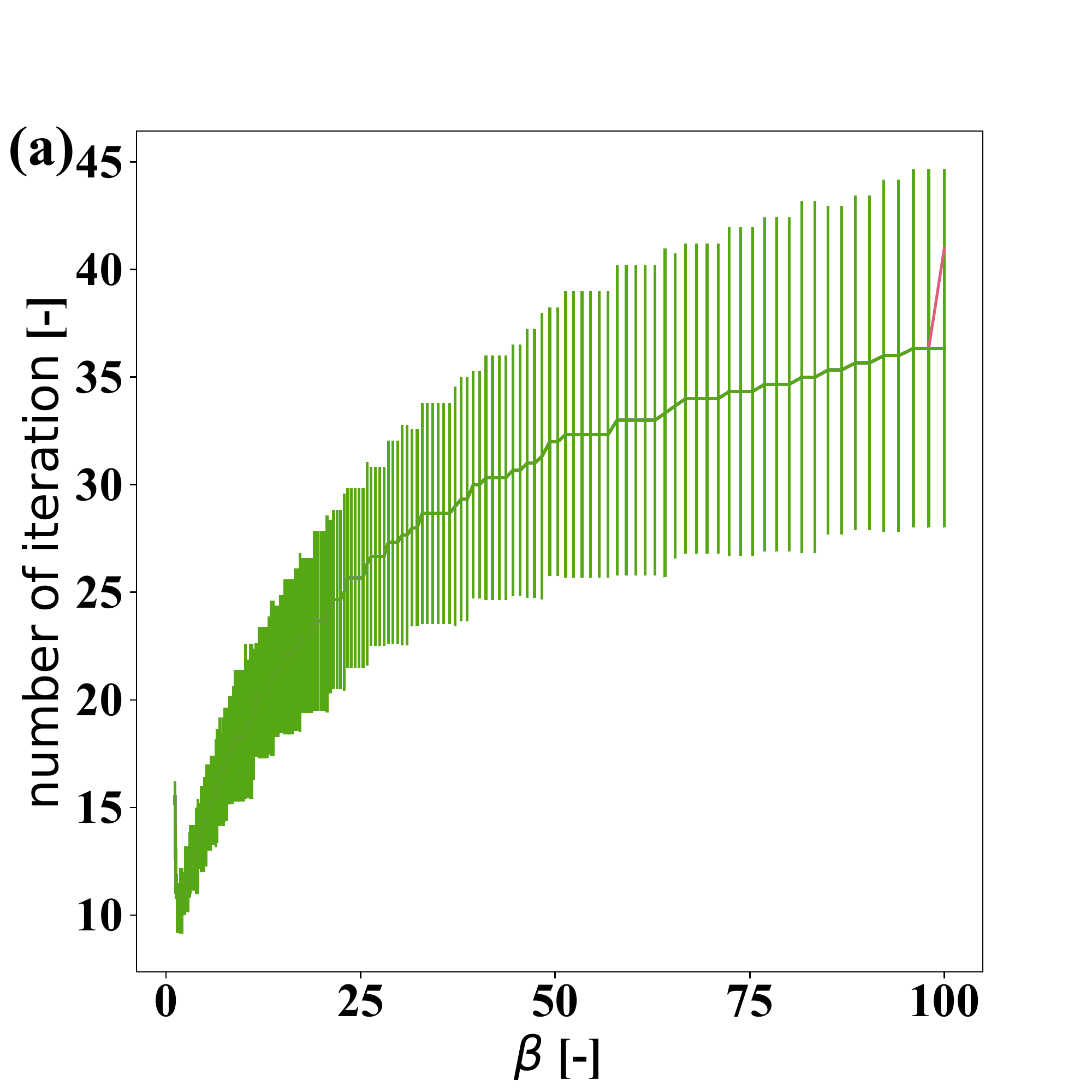}
        \includegraphics[width=0.3\textwidth]{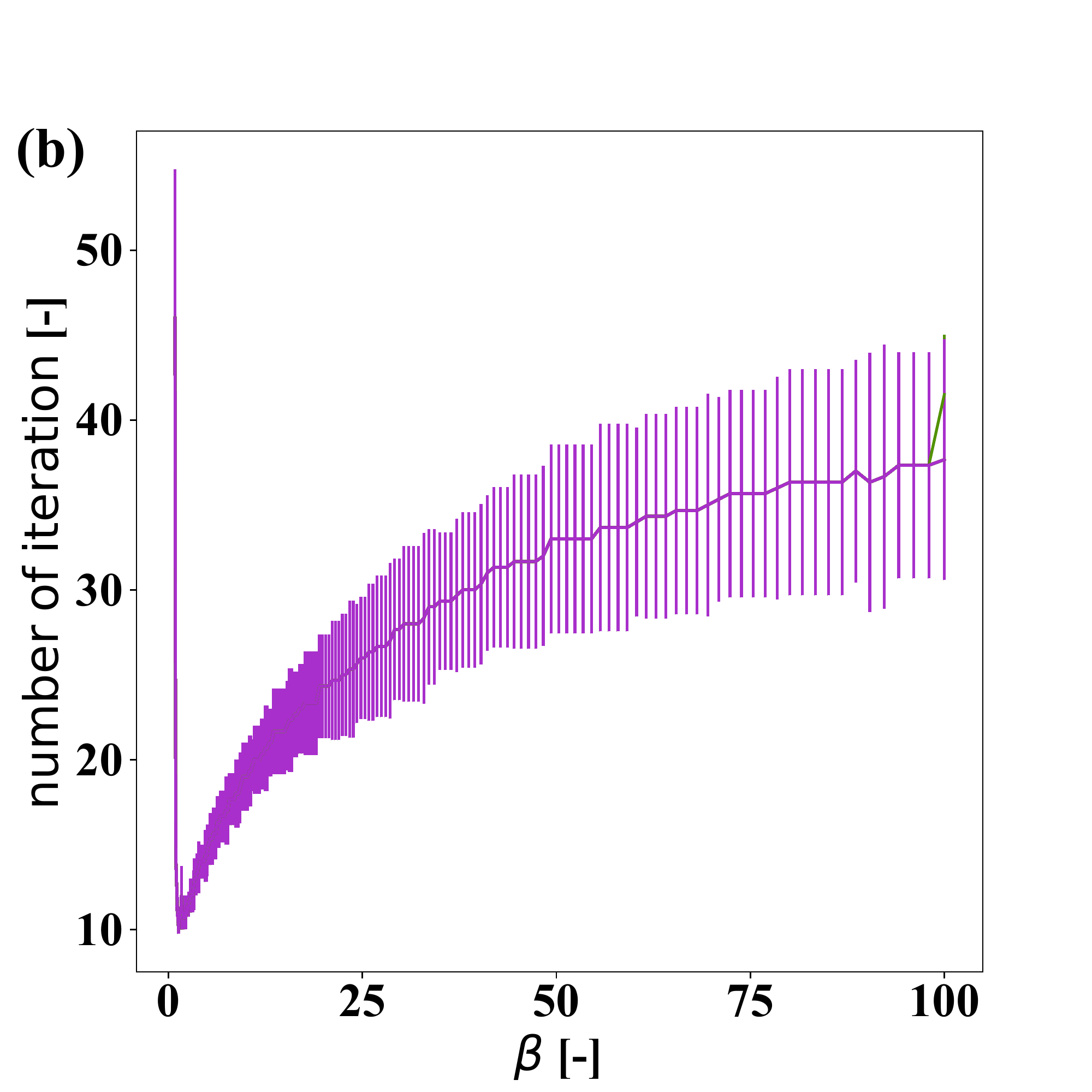}\\
        \includegraphics[width=0.3\textwidth]{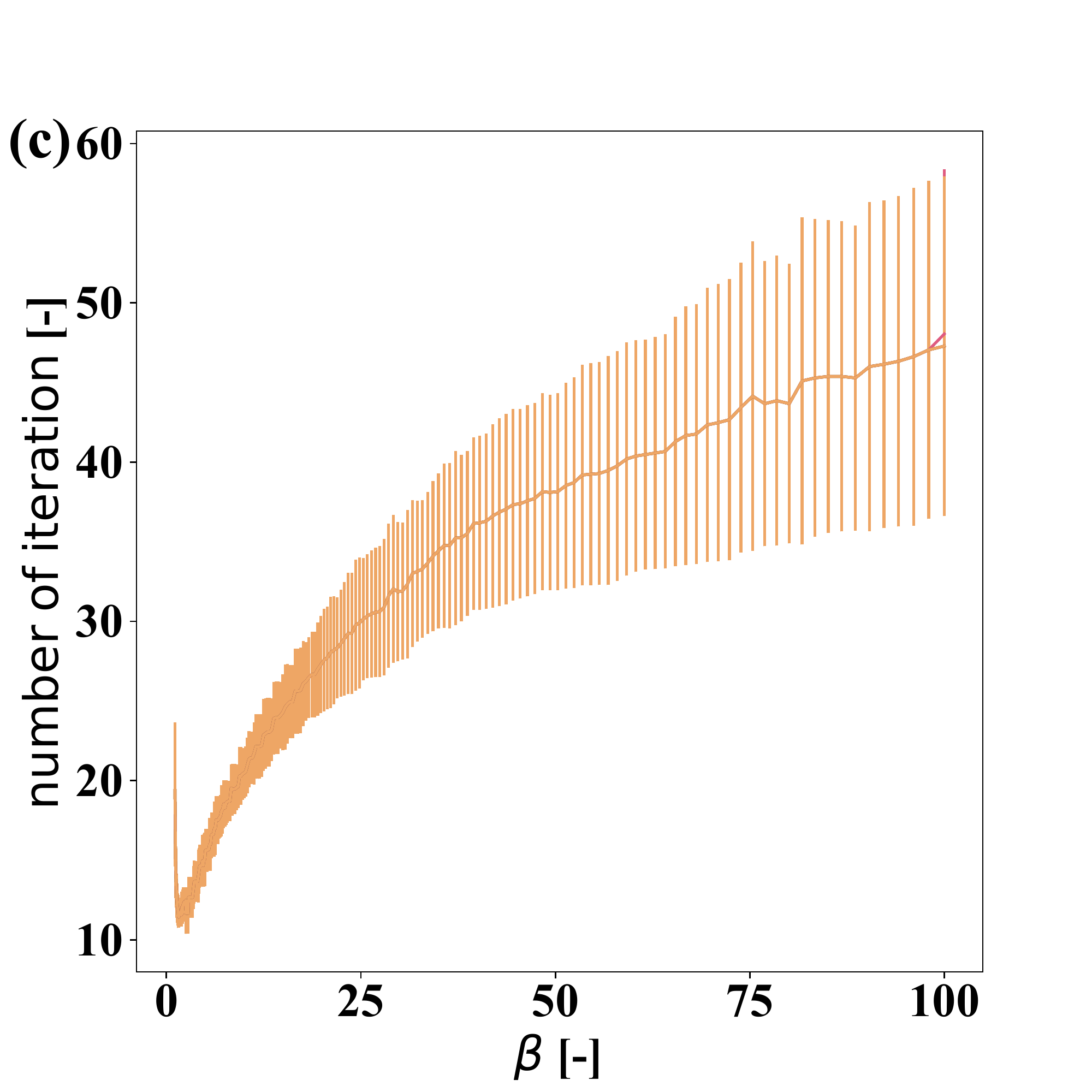}
        \includegraphics[width=0.3\textwidth]{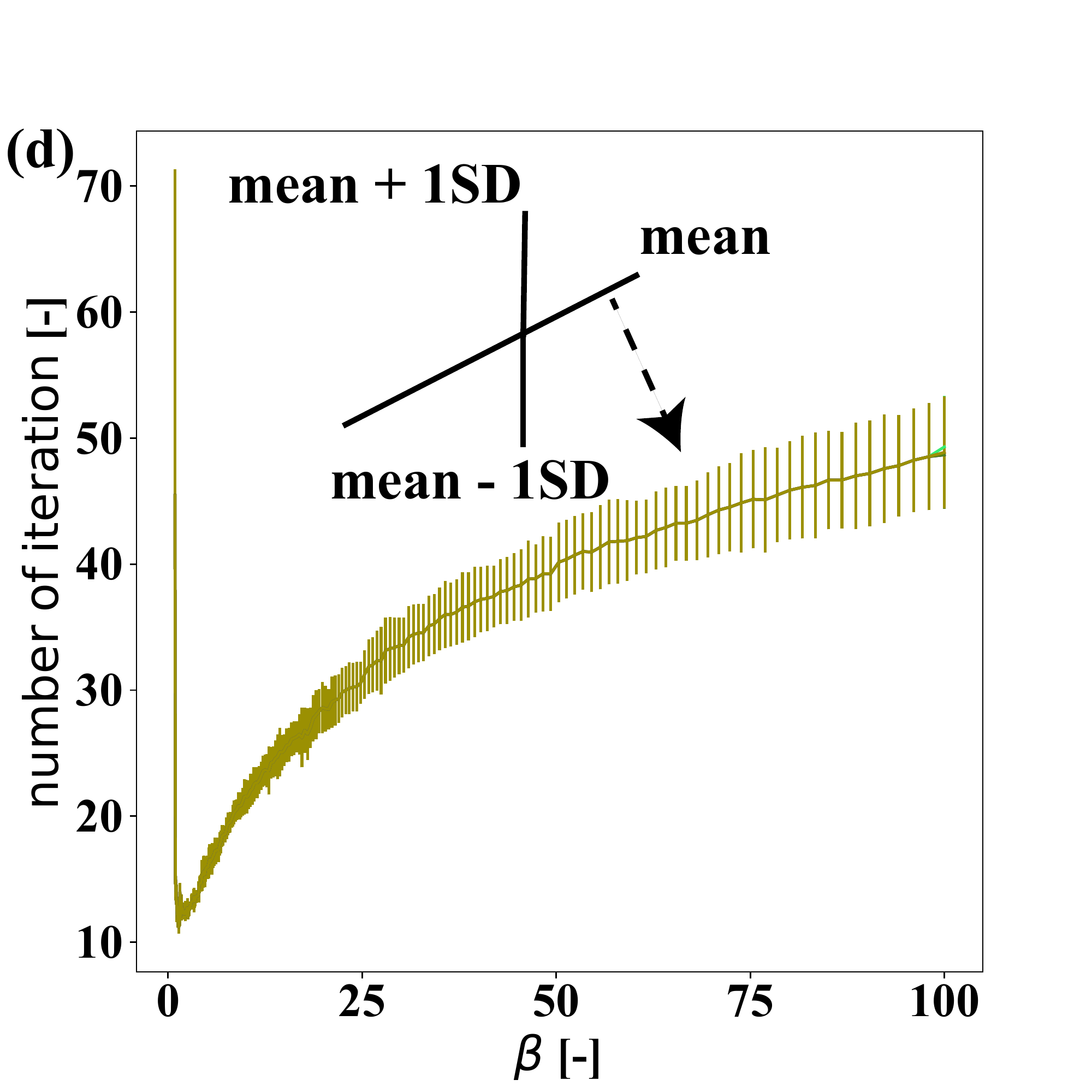}
   \caption{Number of iterations for linear iterative solver with 
   heterogeneous $\boldsymbol{{\kappa}}$. \textbf{(a)} and \textbf{(b)} are the results for the continuous solution \eqref{eq:darcy_exact_solution} by using SIPG and  IIPG, respectivley. \textbf{(c)} and \textbf{(d)} are the results for the discontinuous solution \eqref{eq:dis_exact_solution}  by using SIPG and  IIPG, respectivley. Note that each line represents different value of ${\kappa}$ for the continuous solution, and ${\kappa}_1$ and ${\kappa}_2$ for the discontinuous solution. The error bar shows mean and standard derivation ($\pm$ 1 SD) of each number of iteration see (\textbf{d}) for an explanation.}
   \label{fig:ell_het_noi}
\end{figure}

Similar to the results for the continuous solution presented in the previous section \ref{subsubsec_512}, the choice of $\beta$ influences the number of linear iterative solver significantly, as illustrated in Figure \ref{fig:ell_het_noi}. In short, when $\beta$ is increased, the number of iteration increases, while the number of iteration increases sharply before the solution becomes unstable.

To predict the optimal $\beta$ for the iterative solver, we employ a similar approach that was used for the continuous solution in section \ref{subsubsec_512}. 
First, we evaluate each independent variable using the chi-squared test. 
In this example, we have total 57,835 data points. The result for the chi-squared test  are provided in Table \ref{tab:linear_regression_p_value_dis}.  
It is observed that $\theta$, $\beta$, and $h$ have a p-value of less than 0.025. Hence, we include these variables as independent variables for further predictive model development.
\begin{table}[!ht]
  \centering
    \begin{tabular}{|c|c|}
    \hline
    Variable & p-value \\
    \hline
    $\theta$ & $\approx$ 0.00 \\
    \hline
    $\kappa_0$  & $\approx$ 1.00 \\
    \hline
    $\kappa_1$      & $\approx$ 1.00 \\
    
    \hline
    $\beta$  & $\approx$ 0.00 \\
    \hline
    $h$     & $\approx$ 0.00 \\
    \hline
    
    \end{tabular}%
     \caption{Elliptic equation with discontinuous exact solution: p-value results for each explanatory variable}
  \label{tab:linear_regression_p_value_dis}%
\end{table}%
Then, the developed multi-variable regression \cite{scikit-learn} reads:
\begin{equation} \label{eq:ell_regression_dis}
\begin{split}
\mathrm{number \ of \ iteration} = \alpha 
    + \gamma_{1} \times \theta 
    + \gamma_{2} \times \beta    
    + \gamma_{3} \times h,
\end{split}
\end{equation}
where $\alpha = 19.59$, $\gamma_{1} = -0.86$, $\gamma_{2} = 0.42$, and $\gamma_{3} = -19.51$. Similar to the previous equation \eqref{eq:ell_regression_con}, these parameters provide the minimum value ($\leq 1 \times 10^{-4}$) of MSE value \eqref{eq:loss_mse}.  
Other processes including  the splitting technique and optimization solvers are the same as utilized in the previous model. Finally, the $r^2$ and $\text{EVS}$ obtained for this method is
\begin{equation}
r^2 = 0.84, \ \text{ and } \ \text{EVS} = 0.84
\label{results_linear_2}
\end{equation}

Next, to compare the above results by the nonlinear ANN algorithm, we construct the ANN model using three inputs ($\theta, \beta$, and $h$) and one output (number of iterations) as shown in Figure \ref{fig:ell_ann_fig_dis}. 
The number of hidden layers ($N_{hl}$) and a number of neurons ($N_n$) are used for tuning the hyperparameters. 
ReLU is used as an activation function for each neuron of the hidden layer. 
ADAM and MSE \eqref{eq:loss_mse} are employed for minimization method and loss function, respectively.

\tikzset{%
  every neuron/.style={
    circle,
    draw,
    minimum size=0.1cm
  },
  neuron missing/.style={
    draw=none, 
    scale=1.5,
    text height=0.333cm,
    execute at begin node=\color{black}$\vdots$
  },
}
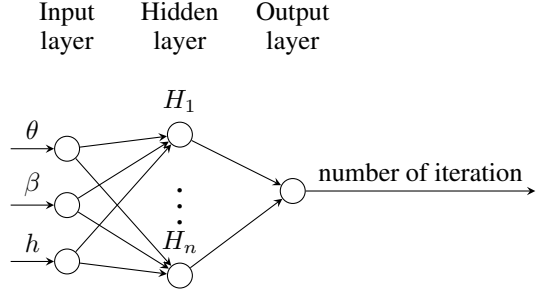
\begin{figure}[!h]
    \centering
    \begin{tikzpicture}[x=1.5cm, y=1.5cm, >=stealth,scale=0.5]
\foreach \m/\l [count=\y] in {1,2,3}
  \node [every neuron/.try, neuron \m/.try] (input-\m) at (0,1.5-\y) {};

\foreach \m [count=\y] in {1,missing,2}
  \node [every neuron/.try, neuron \m/.try ] (hidden-\m) at (2,2-\y*1.25) {};

\foreach \m [count=\y] in {1}
  \node [every neuron/.try, neuron \m/.try ] (output-\m) at (4,0.75-\y) {};

\draw [<-] (input-1) -- ++(-1,0)
    node [above, midway] {$\theta$};
\draw [<-] (input-2) -- ++(-1,0)
    node [above, midway] {$\beta$};
\draw [<-] (input-3) -- ++(-1,0)
    node [above, midway] {$h$};


\foreach \l [count=\i] in {1,n}
  \node [above] at (hidden-\i.north) {$H_\l$};

\draw [->] (output-1) -- ++(4.3,0)
    node [above, midway] {number of iteration};

\foreach \i in {1,...,3}
  \foreach \j in {1,...,2}
    \draw [->] (input-\i) -- (hidden-\j);

\foreach \i in {1,2}
  \foreach \j in {1}
    \draw [->] (hidden-\i) -- (output-\j);

\foreach \l [count=\x from 0] in {Input, Hidden, Output}
  \node [align=center, above] at (\x*2,2) {\l \\ layer};
\end{tikzpicture}
    \caption{Neural network architecture used for the elliptic problem with discontinuous exact solution. The number of hidden layers, $N_{hl}$, and the number of neuron for each hidden layer, $N_n$, are used as the sensitivity analysis parameters. $H_1$ and $H_n$ represent the numbering of each neuron in each hidden layer.}
    \label{fig:ell_ann_fig_dis}
\end{figure}

Table \ref{table:ell_hyper_nl_nn_dis} presents that the MSE of the validation set is decreased as $N_{hl}$ and $N_{n}$ are increased until $N_{hl} = 8$ and $N_{n} = 40$. Hence, we select $N_{hl} = 8$ and $N_{n} = 40$ for the test set. Then we obtain the following final results  
\begin{equation}
r^2 = 0.98, \ \text{ and } \ \text{EVS} = 0.98
\label{results_nonlinear_2}
\end{equation}
The above results from the nonlinear ANN algorithms outperform the linear multi-variable regression \eqref{results_linear_2}. From results of the section \ref{subsec_elliptic}, we can infer that the performance of the nonlinear approximation algorithm is better than the linear one; as a result, the relationship between the number of iteration and its dependent variables is nonlinear.


\begin{table}[!ht]
\centering
\begin{tabular}{|c|c|c|c|c|}
\hline
\backslashbox{$N_{hl}$}{$N_n$}   & 10      & 20      & 40     & 80       \\ \hline

2     & 5.49  & 1.57  & 3.16  & 1.67 \\
    \hline
    4     & 1.93  & 1.53  & 1.71  & 1.63 \\
    \hline
    8     & 1.44  & 2.34  & 1.43  & 1.56 \\
    \hline

\end{tabular}
\caption{Elliptic equation with discontinuous exact solution: Mean squared error of the validation set for different number of hidden layers $N_{hl}$ and different number of neurons per layer $N_{n}$}
\label{table:ell_hyper_nl_nn_dis}
\end{table}

\subsection{Effect and optimal choice of interior penalty parameter for Biot's equations}
\label{subsec_biot}

In this example, we aim to investigate the effect of the interior penalty $\beta$ on the solution quality of the Biot's equations, where the elliptic flow equation is coupled with the solid mechanics as described in Section \ref{subsec_poro}. 
Although, employing DG approximation for the flow equation eliminates any spurious oscillations that are observed when  the continuous Galerkin (CG) is used (especially  at material interfaces where a large conductivity ($\boldsymbol{\kappa}$) contrast is located) as presented in \cite{choo2018enriched,Kadeethum2019},  
the quality of DG solutions may be influenced by the choice of $\beta$.
Thus, we employ the machine learning algorithm to find the optimal choice of $\beta$ to avoid any instabilities upon the given physical and numerical parameters. In the following problems, we compare the performance of the linear logistic regression and nonlinear classification ANN, where the predicted values are binaries.

In the computational domain $\Omega = \left[ 0,1\right]^1$, the geometry and boundary conditions are shown in Figure \ref{fig:biot_2mat}a. 
Here, $\boldsymbol{\kappa}$ is defined as:
\begin{equation}\label{eq:biot_kappa}
\boldsymbol{\kappa} :=
\begin{cases}
{\kappa}_1 \boldsymbol{I} \quad & \text { if }  0.0 \leq x \leq 0.5, \\
{\kappa}_2 \boldsymbol{I} \quad & \text { if }  1.0 \geq x > 0.5, \\
\end{cases} \\
\end{equation}
and we define the ratio between $\kappa_2$ and $\kappa_1$ as
\begin{equation}\label{eq:k_mult}
\kappa_{mult} := \frac{{\kappa}_2 }{{\kappa}_1 }.
\end{equation}

\subsubsection{Effect and optimal choice of interior penalty parameter for Biot's system. }

In this section, we study the optimal choice of interior penalty parameter $\beta$ on the solution quality for the Biot's system.  
The physical parameters are set as 
$\mu=10^{-6}$ $\mathrm{kPa.s}$, 
$\rho=1000$ $\mathrm{kg} / \mathrm{m}^{3}$, 
$K=1000$ $\mathrm{kPa}$, 
$K_{\mathrm{s}} \approx \infty$ $\mathrm{kPa}$, which leads to $\alpha \approx 1$, and $v=0.25$.
In addition, we note that  $\kappa_{1}=10^{-12}$$\mathrm{m}^{2}$ and 
$\kappa_{2}=10^{-16}$$\mathrm{m}^{2}$, and Lam\'e coefficeints $\lambda_l$ and $\mu_l$ are calculated by the following equations:
\begin{equation}
\lambda_{l}=\frac{3 K v}{1+v}, \ \text{ and } \ 
\mu_{l}=\frac{3 K(1-2 v)}{2(1+v)}.
\end{equation}
The numerical parameters are given as $h=0.05$ $m$ and
$\Delta t^n$ $=$ 1.0 sec, and 
the boundary conditions are set to $\boldsymbol{\sigma_{D}}=[0, 1]$ $\mathrm{kPa}$ and $p_{D}=0$ $\mathrm{Pa}$. 
LU direct solver and SIPG ($\theta=1$) is used to solve the discretized system.

\begin{figure}[h!]
   \centering
        \includegraphics[width=0.4\textwidth]{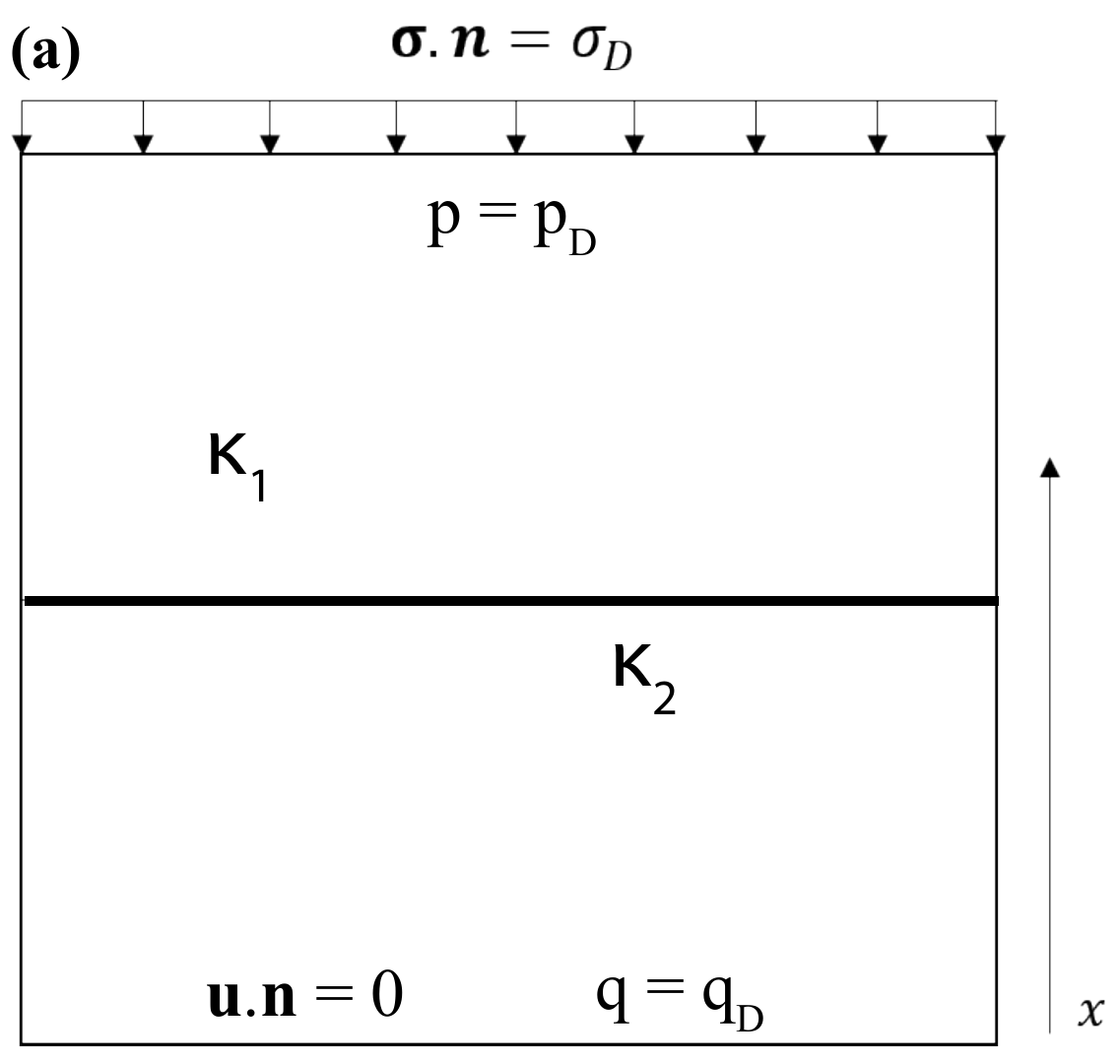}
        \includegraphics[width=0.4\textwidth]{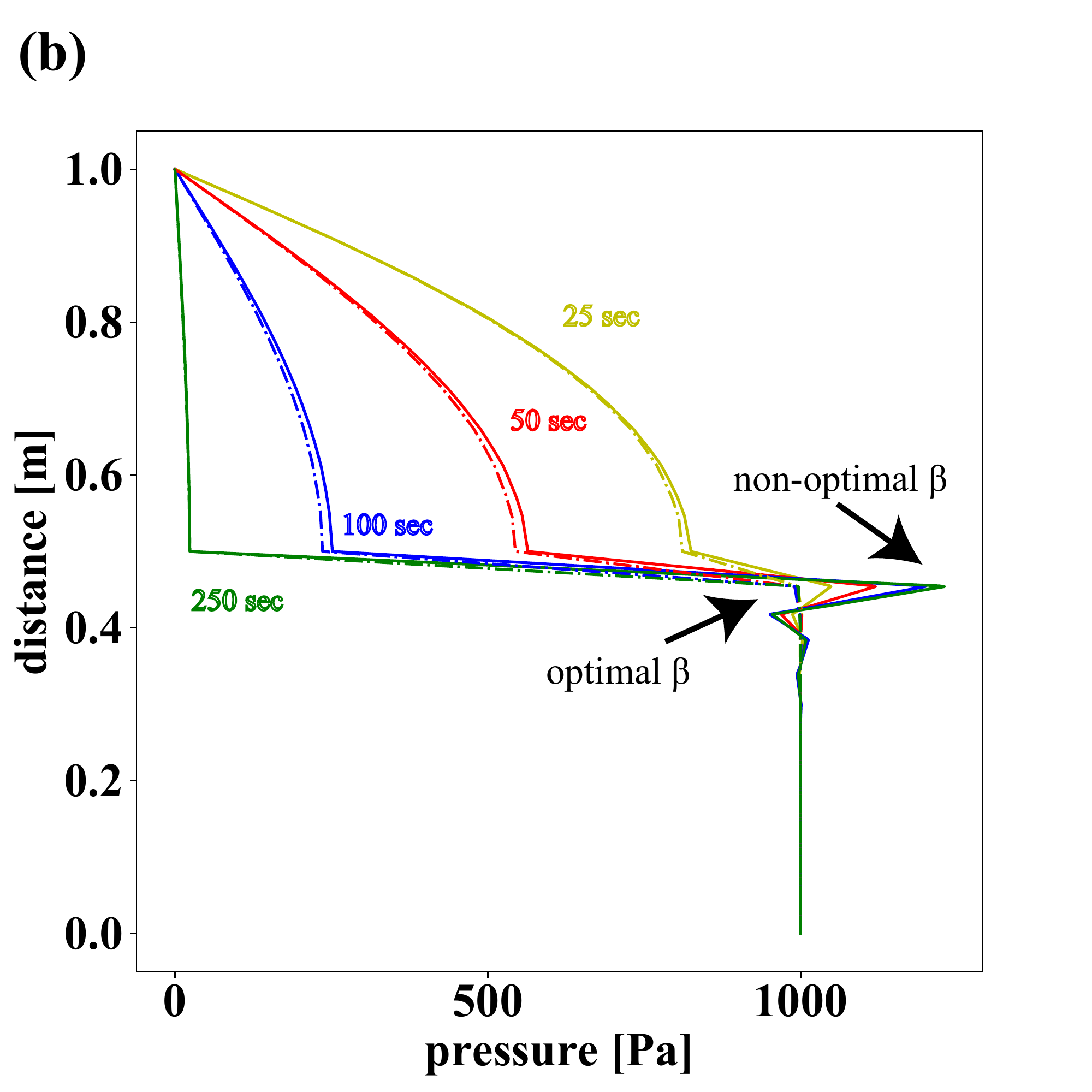}
   \caption{(a) geometry and boundary conditions used in Biot's equations study and (b) pressure results for an example of the effect of $\beta$ on the solution quality}
   \label{fig:biot_2mat}
\end{figure}

For example, the numerical simulation results by comparing 
$\beta = 1.1$ and $\beta = 50.0$ 
are presented in Figure \ref{fig:biot_2mat}b.  
Figure \ref{fig:biot_2mat}b illustrates that the choice of $\beta$ can lead to different qualities of pressure solution, i.e. in case of $\beta = 1.1$ the pressure solution exhibit no spurious pressure oscillations while  the oscillations appear when $\beta = 50.0$.
Note that when $\beta$ is too small, the solution may become also unstable, as illustrated in the previous section for the elliptic problem and discussed in \cite{choo2018enriched} for the Biot's equations.

Next, we investigate to find the optimal choice of $\beta$ using the procedure illustrated in Algorithm \ref{alg:biot_procedure}. 
The ranges of the input values for all the test cases are given as;  ${\kappa}_1$ = $[1.0 \times 10^{-14},1.0 \times 10^{-8}]$, ${\kappa}_2$ = $[1.0 \times 10^{-17},1.0 \times 10^{-14}]$, ${\kappa}_{mult}$ = $[1.0 \times 10^{-8},1.0 \times 10^{1}]$, $\beta$ = $[4.5 \times 10^{-2},2.0 \times 10^{-8}]$, and $h$ = $[7.8 \times 10^{-3},6.25 \times 10^{-2}]$.

To determine the quality of the numerical solution, we define
the two types of solution quality. If the approximated solution is stable and smooth with no spurious pressure oscillations, we denote as  `good.'
If we observe any spurious pressure oscillations from the nonstable approximated solution, we denote as `bad.' Thus, in this case, we utilize the bool type variable {\textsc{Bool\_Quality}} for binary classification which $1$ indicates `good' and  $0$ indicates `bad'.

\begin{algorithm}[!ht]
\caption{Investigation procedure for the Biot's equations}
\label{alg:biot_procedure}
\begin{algorithmic}
\STATE{Initialize sets of each variables; 
{\textsc{Bool\_Quality}}, ${\kappa}$, $\kappa_{mult}$, $h$, and  $\beta$}

\FOR{$i$ $<$ $n_{{\kappa}_1}$, where $n_{{\kappa}_1}$ is the size of the specified ${\kappa}_1$ list}

\STATE{Assign $\boldsymbol{\kappa}_1:={\kappa}_1\left[i\right]\boldsymbol{I}$}

\FOR{$j$ $<$ $n_{\kappa_{mult}}$, where $n_{\kappa_{mult}}$ is the size of the specified $\kappa_{mult}$ list}

\STATE{Assign $\boldsymbol{\kappa}_2:=\kappa_{mult}\left[j\right] \times \boldsymbol{\kappa}_1$}

\FOR{$k$ $<$ $n_h$, where $n_h$ is the size of the specified $h$ list}

\STATE{Assign $h:=h\left[k\right]$}

\FOR{$l$ $<$ $n_{\beta}$, where $n_{\beta}$ is the size of the specified $\beta$ list}

\STATE{Assign $\beta:=\beta\left[l\right]$}
\STATE{Solve the coupled Biot's system: \eqref{eq:dgscheme} and \eqref{eq:cgscheme}.}

\IF{linear solver converges} 
\IF{nonphysical spurious oscillation is detected} 
\STATE{\textsc{Bool\_Quality} = 0}
\ELSE 
\STATE{\textsc{Bool\_Quality} = 1}
\ENDIF
\ELSE 
\STATE{\textsc{Bool\_Quality} = 0}
\ENDIF
\ENDFOR
\ENDFOR
\ENDFOR
\ENDFOR
\end{algorithmic}
\end{algorithm}

Similar to the previous sections, we begin with the chi-squared test to find the statistically significant explanatory variables. In total, we have 14,141 cases (data points) with 3,927 `good'({\textsc{Bool\_Quality} = 1}) solutions and 10,214 `bad'({\textsc{Bool\_Quality} = 0}) solutions. 
The chi-squared test result is presented in Table \ref{tab:linear_regression_p_value_biot}, and it shows that all variables, ${\kappa}_1$, ${\kappa}_2$, $\kappa_{mult}$, and $h$, have p-value less than 0.025. Therefore, these variables are included as independent variables to develop the following predictive models.
\begin{table}[!h]
  \centering
    \begin{tabular}{|c|c|}
    \hline
    Variable & p-value \\
    \hline
    ${\kappa}_1$ & $\approx$ 0.00 \\
    \hline
    ${\kappa}_2$ & $\approx$ 0.00 \\
    \hline
    $\kappa_{mult}$ & $\approx$ 0.00 \\
    \hline
    $\beta$  & $\approx$ 0.00 \\
    \hline
    $h$     & $\approx$ 0.01 \\
    \hline
    \end{tabular}
\caption{Biot's equations: p-value results for each explanatory variables.}
\label{tab:linear_regression_p_value_biot}%
\end{table}%
As discussed in both sections \ref{sec:alg_1} and \ref{sec:alg_2}, the data sets are split by  training, validation, and test sets  
using the splitting ratio 
$\left[ 0.8, 0.1 ,0.1\right]$. Thus, the number of  training sets, validation sets, and test sets are  $0.8 \times 14,141$, $0.1 \times 14,141$, and  $0.1 \times 14,141$, respectively.
We use the training set to train the linear and nonlinear machine learning algorithms. The validation set is for tuning the hyperparameters for the nonlinear ANN models, and the test set is for comparing performances between the linear and nonlinear algorithms.


First, the multi-variable logistic regression \cite{scikit-learn} is defined as:
\begin{equation} \label{eq:biot_regression_dis}
\begin{split}
\log\left(\frac{\mathbb{P}(\mathbb{O}=1)}{1-(\mathbb{P}(\mathbb{O}=1))}\right) = \alpha 
    + \gamma_{1} \times \boldsymbol{\kappa}_1 
    + \gamma_{2} \times \boldsymbol{\kappa}_2 
    + \gamma_{3} \times \kappa_{mult} 
    + \gamma_{4} \times \beta    
    + \gamma_{5} \times h,
\end{split}
\end{equation}
where
\begin{equation}\label{eq:exact_solution}
\mathbb{O}=
\begin{cases}
1, & \text{if} \ \mathbb{P}(\mathbb{O}=1) \ge 0.5 \\
0, & \text{if} \ \mathbb{P}(\mathbb{O}=1) < 0.5.
\end{cases}
\end{equation}
Here, $\alpha = -1.19$, $\gamma_{1} = -0.06$, $\gamma_{2} = 0.68$, $\gamma_{3} = 5.55$, $\gamma_{4} = -5.49$, and $\gamma_{5} = 0.32$. These parameters provide the minimum value ($\leq 1 \times 10^{-4}$) of BCE value \eqref{eq:loss_bce}. 

After applying the algorithm explained in section
\ref{sec:alg} and \ref{sec:alg_1}, the computed accuracy \eqref{eq:loss_acc} of the logistic regression model is 
$$
\text{ACC} = 0.80,
$$
and the confusion matrix is presented in Table \ref{tab:biot_cm_logistic}. 
In table \ref{tab:biot_cm_logistic}, we observe that the number of 
`false positive' is much higher than to that of `false negative,' which may result in the bad solution obtained from the finite element model. To elaborate, when our model creates 'false positive,' it means that we expect the simulation results to be stable and contain no oscillation; however, in fact, the solution quality is bad.

\begin{table}[!ht]
  \centering
    \begin{tabular}{|c|c|c|c|}
    \hline
    \multicolumn{2}{|c|}{\multirow{2}[4]{*}{total test set = 1415}} & \multicolumn{2}{c|}{Test set values} \\
\cline{3-4}    \multicolumn{2}{|c|}{} & \multicolumn{1}{l|}{Good (1)} & \multicolumn{1}{l|}{Bad (0)} \\
    \hline
    \multirow{2}[4]{*}{Predicted values} & Good (1) & 172   & 210 \\
\cline{2-4}          & Bad (0) & 74   & 959 \\
    \hline
    \end{tabular}%
    \caption{Biot's equations: Confusion matrix of the logistic regression for the test set}
  \label{tab:biot_cm_logistic}%
\end{table}
\begin{table}[!h]
\centering
\begin{tabular}{|c|c|c|c|c|c|}
\hline
\backslashbox{$N_{hl}$}{$N_n$}   & 10      & 20      & 40     & 80   & 120     \\ \hline

    2     & 0.89  & 0.93  & 0.93  & 0.93  & 0.92 \\
    \hline
    4     & 0.89  & 0.93  & 0.93  & 0.93  & 0.93 \\
    \hline
    8     & 0.88  & 0.92  & 0.93  & 0.93  & 0.93 \\
    \hline
    16    & 0.73  & 0.73  & 0.73  & 0.27  & 0.27 \\
    \hline
    32    & 0.73  & 0.73  & 0.73  & 0.73  & 0.73 \\
    \hline
\end{tabular}
\caption{Biot's equations: Accuracy of the validation set for different number of hidden layers $N_{hl}$ and different number of neurons per layer $N_{n}$}
\label{table:biot_hyper_nl_nn_dis}
\end{table}
Secondly, we develop the classification ANN model utilizing five input ( $\boldsymbol{\kappa}_1, \boldsymbol{\kappa}_2, \kappa_{mult}, \beta,$ and $h$)  and one output ({\textsc{Bool\_Quality}})as shown in Figure \ref{fig:biot_ann_fig} for this problem. Table \ref{table:biot_hyper_nl_nn_dis} illustrates the result for hyperparameters tuning, and it illustrates that the ANN predictive performance is improved as $N_{hl}$ and $N_{n}$ are increased up until $N_{hl} = 4$ and $N_{n} = 80$. Hence, we use $N_{hl} = 4$ and $N_{n} = 80$ to set  the hyperparameters, and we compare the nonlinear classification ANN and the logistic regression models' performance. Here, ReLU is used as an activation function for each neuron of the hidden layer, and the Sigmoid activation function is used for the output layer. ADAM and $\text{BCE}$ \eqref{eq:loss_bce} are employed for minimization method and loss function, respectively.

\tikzset{%
  every neuron/.style={
    circle,
    draw,
    minimum size=0.1cm
  },
  neuron missing/.style={
    draw=none, 
    scale=1.5,
    text height=0.333cm,
    execute at begin node=\color{black}$\vdots$
  },
}
\begin{figure}
    \centering
    \begin{tikzpicture}[x=1.5cm, y=1.5cm, >=stealth, scale=0.5]
\foreach \m/\l [count=\y] in {1,2,3,4,5}
  \node [every neuron/.try, neuron \m/.try] (input-\m) at (0,2.5-\y) {};

\foreach \m [count=\y] in {1,missing,2}
  \node [every neuron/.try, neuron \m/.try ] (hidden-\m) at (2,2-\y*1.25) {};

\foreach \m [count=\y] in {1}
  \node [every neuron/.try, neuron \m/.try ] (output-\m) at (4,0.75-\y) {};

\draw [<-] (input-1) -- ++(-1,0)
    node [above, midway] {$\boldsymbol{\kappa}_1$};
\draw [<-] (input-2) -- ++(-1,0)
    node [above, midway] {$\boldsymbol{\kappa}_2$};
\draw [<-] (input-3) -- ++(-1,0)
    node [above, midway] {$\kappa_{mult}$};
\draw [<-] (input-4) -- ++(-1,0)
    node [above, midway] {$\beta$};
\draw [<-] (input-5) -- ++(-1,0)
    node [above, midway] {$h$};


\foreach \l [count=\i] in {1,n}
  \node [above] at (hidden-\i.north) {$H_\l$};

\draw [->] (output-1) -- ++(4.2,0)
    node [above, midway] {\textsc{Bool\_Quality}};

\foreach \i in {1,...,5}
  \foreach \j in {1,...,2}
    \draw [->] (input-\i) -- (hidden-\j);

\foreach \i in {1,2}
  \foreach \j in {1}
    \draw [->] (hidden-\i) -- (output-\j);

\foreach \l [count=\x from 0] in {Input, Hidden, Ouput}
  \node [align=center, above] at (\x*2,2) {\l \\ layer};
\end{tikzpicture}
    \caption{Neural network architecture used for the Biot's equation. The number of hidden layers, $N_{hl}$, and the number of neuron for each hidden layer, $N_n$, are used as the sensitivity analysis parameters. $H_1$ and $H_n$ represent the numbering of each neuron in each hidden layer.}
    \label{fig:biot_ann_fig}
\end{figure}
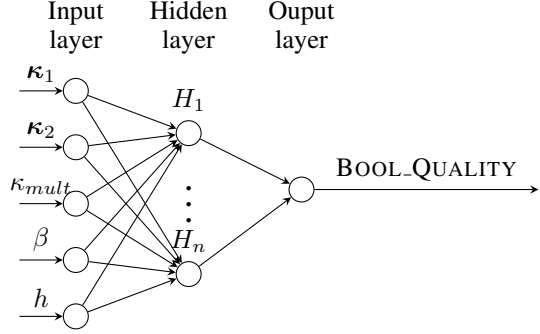

Finally, the computed accuracy value \eqref{eq:loss_acc} of the nonlinear classification ANN using the above test set is 
$$ 
\text{ACC} = 0.93.
$$ 
We note that this value is much higher than that of the logistic regression model. 

Furthermore, the number of `false positive' cases presented in Table \ref{tab:biot_cm_ann} is much lower compared to that of the linear logistic regression algorithm. This characteristic helps to prevent the finite element model from producing bad quality simulation results as discussed previously.

\begin{table}[!ht]
  \centering
    \begin{tabular}{|c|c|c|c|}
    \hline
    \multicolumn{2}{|c|}{\multirow{2}[4]{*}{total test set = 1415}} & \multicolumn{2}{c|}{Test set values} \\
\cline{3-4}    \multicolumn{2}{|c|}{} & \multicolumn{1}{l|}{Good (1)} & \multicolumn{1}{l|}{Bad (0)} \\
    \hline
    \multirow{2}[4]{*}{Predicted values} & Good (1) & 367   & 15 \\
\cline{2-4}          & Bad (0) & 72    & 961 \\
    \hline
    \end{tabular}%
    \caption{Biot's equations: Confusion matrix of the artificial neural network (ANN) for the test set}
  \label{tab:biot_cm_ann}%
\end{table}%


\section{Conclusions}
\label{sec:conclusions}

This paper presents the effect of the choice of the interior penalty parameter of the discontinuous Galerkin finite element methods for the elliptic problems and the Biot’s systems. The optimal choice of the interior penalty parameter results in stable solutions, optimum error convergence rate, less number of iteration for the iterative solver, and eliminates any spurious numerical oscillation in the approximated solutions. We propose the nonlinear approximation algorithms, regression and classification, to predict the optimal choice of the interior penalty parameter. These nonlinear approximation algorithms outperform the classic linear approximation algorithms. Our proposed framework can be beneficial to sensitivity analysis, uncertainty quantification, or data assimilation modelling where many simulations have to be performed with different settings, e.g., mesh size, material properties, or different interior penalty schemes. Moreover, it can be extended to any multiscale multiphysics problems.

\section*{Acknowledgments}
SL is supported by National Science Foundation under Grant No. NSF DMS-1913016.
TK and HM have received funding from the Danish Hydrocarbon Research and Technology Centre under the Advanced Water Flooding program. 
We acknowledge developers and contributors of TensorFlow \cite{tensorflow2015-whitepaper}, Keras \cite{chollet2015keras}, Scikit-learn \cite{scikit-learn}, FEniCS \cite{AlnaesBlechta2015a}, and Multiphenics \cite{Ballarin2019} libraries.

\bibliography{lit}
\end{document}